\documentclass[notitlepage, 11pt]{article}

\usepackage{color}

\newtheorem{theorem}{Theorem}

\newtheorem{definition}{Definition}

\newtheorem{condition}{Condition}

\newtheorem{remark}{Remark}

\newcommand{\noi}{\noindent}
\newcommand{\spa}{\vspace{.2in}}

\newcommand{\E}{\mathbb{E}}
\newcommand{\R}{\mathbb{R}}

\newcommand{\p}{\mathbb{P}}
\newcommand{\N}{\mathbb{N}}
\newcommand{\I}{\mathcal{I}}
\newcommand{\Ir}{\mathbb{I}}

\newcommand{\la}{\lambda}
\newcommand{\sig}{\sigma}

\newcommand{\eps}{\varepsilon}
\newcommand{\e}{\mathbf{e}}
\newcommand{\m}{\mathbf{d}}
\newcommand{\bzeta}{\boldsymbol{\zeta}}

\newcommand{\ph}{\varphi}

\newcommand{\al}{\alpha}

\newcommand{\gam}{\gamma}
\newcommand{\kap}{\kappa}

\newcommand{\del}{\delta}
\newcommand{\om}{\omega}

\newcommand{\Gam}{\mathnormal{\Gamma}}
\newcommand{\Del}{\mathnormal{\Delta}}

\newcommand{\La}{\mathnormal{\Lambda}}

\newcommand{\Ups}{\mathnormal{\Upsilon}}
\newcommand{\Om}{\mathnormal{\Omega}}

\newcommand{\calI}{{\cal I}}

\newcommand{\calQ}{{\cal Q}}

\newcommand{\calS}{{\cal S}}
\newcommand{\calT}{{\cal T}}

\newcommand{\bI}{{\mathbf I}}

\newcommand{\frB}{\mathfrak{B}}

\newcommand{\oo}{\overline}

\newcommand{\w}{\wedge}

\newcommand{\To}{\Rightarrow}

\newcommand{\iy}{\infty}
\newcommand{\be}{\begin{equation}}
\newcommand{\ee}{\end{equation}}
\newcommand{\osc}{\text{osc}}
\newcommand{\ds}{\displaystyle}
\newcommand{\A}{{\cal A}}

\newcommand{\Ri}{\R^{2\bI}}

\usepackage{latexsym}
\usepackage{amsmath, amssymb, amsfonts}

\numberwithin{equation}{section}

\title{Risk-sensitive control for the multi-class many server queue in the moderate deviation regime}
\author{ Anup Biswas\\ \\ \\
Department of Electrical Engineering\\
Technion--Israel Institute of Technology\\
Haifa 32000, Israel\\
email:{\it anup@ee.technion.ac.il}}

\date{}

\begin{document}

\maketitle

\begin{abstract}
A G/M/N queue is considered in the moderate deviation heavy traffic regime. The rate function for the customers-in-system
process is obtained for the single class model. A risk-sensitive type control problem is considered for multi-class G/M/N
model under the moderate deviation scaling and shown that the optimal control problem is related to a differential game problem.

\spa

\noi{\bf AMS subject classifications:}\,\, 60F10, 60K25, 49N70, 93E20

\spa

\noi{\bf Keywords:}\,\,
Risk-sensitive control, large deviations, moderate deviations, differential games,
multi-class many-server queue, G/M/N, heavy traffic

\end{abstract}

\section{Introduction}
Studying scaling limit is an established tradition in queuing theory. These include heavy traffic approximation that depends on central limit theorem (CLT) and large deviation (LD) approximations. Another interesting scaling, considered in queuing network, is moderate
deviation (MD) scaling which includes an intermediate scaling of CLT and LD. MD scaling are consider when the queuing network is critically loaded. Therefore MD can be seen as a LD analogue for heavy traffic set up. Also some control problems in MD regime
have interesting characteristic that also appears in the asymptotic regime alluded to LD and heavy traffic approximations \cite{Atar-Bis}.

There have been several works on MD scaling without dynamic control aspect. LD and MD for renewal processes are proved in \cite{Puhal-Whit}. Later in \cite{Puhal-1999}, Puhalskii obtains the MD principle for the queue length and waiting-time processes for
a single class single server network. Majewski \cite{Majewski06} considers feedforward multi-class network with priority and obtains
the MD asymptotics for waiting time, idle time, queue length, departure and sojourn time processes. We refer to \cite{Ganesh, Wischik} for various interesting aspects of MD regime. A dynamic control problem for multi-class G/G/1 queue in the MD regime
is considered in \cite{Atar-Bis} where the authors point out some interesting features of the problem similar to other
asymptotic regimes.

So far MD asymptotics have not been considered in many server queuing network. In this article we introduce the MD principle
for the customers-in-system process in a many server network. We consider a single class G/M/N queuing network where the arrivals are given by a general renewal process and the service requirements of the customers are exponentially distributed. We show that the rate function for customers-in-system process in the MD regime changes depending on the growth rate of number of servers $N$ compare to the arrival rate $\la^n$. It is shown that if $N=o(\la^n)$ then the rate function for the customers-in-system process in the MD regime  is governed by a Skorohod map. But if $\frac{N}{\la^n}\not\to 0$ as $n\to\iy$, the governing dynamics for the rate function are not reflection maps. It is worthwhile to mention that this problem can be seen as the MD analogue of the scaling considered by Halfin and Whitt for G/M/N queuing network in \cite{halfin-whitt}.
One may wish to consider the MD analysis for G/G/N queuing network
but the problem is harder as one needs to consider an infinite dimensional set up for the problem.

We also consider a risk-sensitive type control problem for a multi-class G/M/N network when $N=o(\la^n)$. We consider $\bI$ different
customer classes arriving to a parallel server system following $\bI$ independent renewal processes. Service time distributions are
exponential with class dependent parameters. Each customer is served by one of the servers and servers are not allowed to serve
more than one customers at the same time. The problem is to control $B^n=(B^n_1,\ldots,B^n_\bI)$ where $B^n_i$ denotes the number of class$-i$ customers receiving service, so that the cost is minimized. Denoting by $X^n_i$, the number of class$-i$ jobs in the $n$-th system,
the scaled version is given by $\tilde X^n_i=\frac{X^n_i-\rho_iN}{b_n\sqrt n}$ where $\rho_i$ denotes the limiting traffic intensity
for class$-i$ and $\lim b_n=\iy,\ \lim\frac{b_n}{\sqrt n}=0$. The cost is given by
$$\frac{1}{b_n^2}\E[e^{b^2_n(\int_0^Th(\tilde X^n(s)ds+g(\tilde X^n(T)))}], $$
where $T>0,$ and $h,\ g$ are given nonnegative functions. The $n-$th value function is defined to be the infimum of the above cost where the infimum
is taken over all admissible controls. The goal is to study the limit of these value functions as $n\to\iy$.
This risk-sensitive type of cost  has been studied in literature
for its own importance (see \cite{ata-dup-shw, Atar-Gos-Shw, whittle-book}). One of the important aspect of the exponential cost is that it penalizes the large quantities heavily. This is one of the reason for considering exponential cost attached to the queue length
or customers-in-system processes. Another interesting aspect of working in MD regime is that the limiting differential game (DG) is solvable \cite{Atar-Bis}.

It is also interesting to compare the control problem above with the existing similar control problems (\cite{Atar-Bis, Atar-Gos-Shw}). In \cite{Atar-Gos-Shw}, the authors consider a similar problem (with bounded $h$) for multi-class M/M/N network
in the LD regime. The convergence of the value functions, corresponding to the above optimal control problem, is proved using the Markov property and various martingale estimates on the underlying dynamics.
In \cite{Atar-Bis}, a similar problem is treated for multi-class G/G/1 network and the convergence result is obtained by constructing a particular policy. In both the problems, the servers are allowed to serve more than one customer simultaneously. Since the arrival process here is
given by a general renewal process the underlying state dynamics $X^n$ is not Markov and hence \cite{Atar-Gos-Shw} does not apply.
First of all, our proof technique here does not use any PDE analysis as in \cite{Atar-Gos-Shw}.
Also we do not allow processor sharing. So the set of controls considered in this paper is smaller
than those that are considered in earlier works. The proof of the convergence of the value function for the optimal control problem is divided
into two parts. We first prove the lower bound estimate following similar technique as \cite{Atar-Bis}. The proof for the upper bound
is based on the construction of a particular policy such that the lower bound is asymptotically attained. The construction of
this policy is complicated than that appear in \cite{Atar-Bis} and can be used to improve the control set used in \cite{Atar-Bis}.
We also obtain a simple control that is asymptotically optimal when the cost functions are linear and $N=o(b_n\sqrt n)$. \cite{atagur, atasol} deal with a multi-class G/M/N network under diffusion scaling where $N\approx \sqrt n$.  Our problem can also be thought of as a generalization to these works in risk-sensitive set up. Let us also mention a related work \cite{atar-avi-reiman} where a multi-class scheduling problem is considered
under diffusion scaling.

To summarize the main contribution of the paper, we have (a) introduced the moderate deviation scaling for the many server queues
in heavy traffic regime, (b) shown the convergence of value function for the optimal control problem to a value function of DG,
(c) considered a smaller class of {\it admissible} control which can also be used to improve
the results in \cite{Atar-Bis}, (d) given a simple policy when the cost functions are linear and $N=o(b_n\sqrt n)$.

\spa

\noi{\bf Notations:}\
For a positive integer $k$ and $a,b\in\R^k$, $a\cdot b$ denotes the usual
scalar product, while $\|\cdot\|$ denotes the Euclidean norm. The interval $[0, \iy)$ is denoted by $\R_+$. For $a\in\R_+$, $\lfloor a\rfloor$ denote the largest integer
less than or equal to $a$. Given $a, b\in\R$, the maximum (minimum) is denoted by $a\vee b$ ($a\wedge b$). We use $a^+\ (a^-)$ for $a\vee 0\ (-a\vee 0)$. Given two sequences $\{a_n\},\ \{b_n\}$, $a_n=o(b_n)$ means $\limsup\frac{a_n}{b_n}=0$.
By $\R^k_+$ we denote the nonnegative orthant of the Euclidean space $\R^k$.
For $T>0$ and a function $f:[0,T]\to\R^k$,
we define $\|f\|^*_t=\sup_{s\in[0,t]}\|f(s)\|$, $t\in[0,T]$. When $k=1$, we write $|f|^*_t$ for $\|f\|^*_t$
and $\|f\|^*$ for $\|f\|^*_T$. $\e(\cdot)$ is used to denote the identity function on $\R$.
For $T\leq\iy$, denote by $C([0,T],\R^k)$ and $D([0,T], \R^k)$ the
spaces of continuous functions $[0,T]\to\R^k$ and respectively, functions that are
right-continuous with finite left limits (RCLL). For fix $T>0$, endow the space $D([0,T],\R^k)$ with the
Skorohod-Prohorov-Lindvall metric or $J_1$ metric,
defined as
$$
\m(\ph,\ph^\prime)=\inf_{f\in\Ups}\Big(\|f\|^\circ\vee\sup_{[0, T]}\|\ph(t)-\ph^\prime(f(t))\| \Big),
\qquad \ph, \ph^\prime\in D([0, T], \R^k),
$$
where $\Ups$ is the set of strictly increasing, continuous functions from $[0, T]$ onto itself,
and
$$\| f\|^\circ=\sup_{0\leq s<t\leq T}\Big|\log\frac{f(t)-f(s)}{t-s}\Big|.$$
As is well known \cite{Bill}, $D([0, T], \R^k)$ is a Polish space under the induced topology. Through out this article, we fix a complete probability space $(\Om, \mathcal{F}, \p)$. All the stochastic processes introduced in this article are defined on $(\Om, \mathcal{F}, \p)$.

The paper is organized as follows. The next paragraph introduces some preliminaries that will be used in this paper. Section \ref{sec2} introduces the moderate deviation principle for the single class G/M/N queue. Section \ref{sec3} is devoted to the study of the multi-class G/M/N queues and the dynamic control problem. Section \ref{sec3.1} introduces the associated dynamic games and states the main results. The proof of the main results are given in Section \ref{sec3.2}. Finally, in Section \ref{sec3.3} we prescribe a simple control
which is asymptotically optimal when the cost functions are linear.

\spa

\noi{\bf Preliminaries:}\
Now we state the definition and properties of large deviation principle (LDP) and Skorohod problem
that will be used in this paper. Given a metric space $\calS$,
a function $\Ir$, defined on $\calS$, is said to ba a  {\it rate function} if the set $\{x\in\calS :\ \Ir(x)\leq a\}$ is compact
for all $a\geq 0,$ and there is a sequence $\{\p_n\}_{n\geq 1}$ of probability measure on the Borel $\sig$-field of $\calS$
(or sequence of random variable $\{X^n\}$ with law $\p_n$) satisfying large deviation principle (LDP) with parameter $a_n\to\infty$ and
rate function $\Ir$ i.e.,
$$ \limsup\frac{1}{a_n}\log\p_n(F)\leq -\inf_{x\in F}\Ir(x),$$
for all closed set $F\subset\calS$, and
$$ \liminf\frac{1}{a_n}\log\p_n(G)\geq -\inf_{x\in G}\Ir(x),$$
for all open set $G\subset\calS$.

One standard way to get new LDP's from an existing one is through {\it contraction mapping principle} which states that
if $\{X^n\}$ satisfies LDP with rate function $\Ir$ and $f$ is a continuous function on $\calS$, then $f(X^n)$ satisfies
LDP with rate function
\be\label{rate}
\Ir_f(y)=\inf_{x:\ y=f(x)}\Ir(x).
\ee
There are several extension to this contraction mapping principle. We refer to \cite{Garcia} for a survey on contraction
mapping principles. In this article, we use an {\it extended contraction mapping principle} which states that if $\{X^n\}$
obeys LDP with rate function $\Ir$, $\{f^n\}$ is a sequence of measurable functions, and if there is a measurable function $f$,
continuous when restricted to the set $\{x: \ \Ir(x)\leq a\},\ a\geq 0$, and $f^n(x^n)\to f(x)$ as $n\to\iy$ whenever $x^n\to x$ and $\Ir(x)<\iy$, then $\{f^n(X^n)\}_{n\geq 1}$ obeys LDP with rate function given by \eqref{rate}.

Our goal in this paper is to study asymptotics of certain value functions and to show that they lead to the value function
of certain differential game problem. This differential game problem is solvable.
In order to define the solution to the game we need do define {\it Skorohod problem}.
\begin{definition}
Let $\psi\in D([0,\iy), \R)$ with $\psi(0)\in R_+$ be given. Then $(\phi^1, \phi^2)$ solves the Skorohod problem for
the data $\psi$ if $\psi(0)=\phi^1(0)$, and for all $t\in[0,\iy)$
\begin{enumerate}
\item $\phi^1(t)=\psi(t)+\phi^2(t)$,
\item $\phi^1(t)\in\R_+$,
\item $\phi^2$ is  nondecreasing,
\item $\int_0^\iy \phi^1(s)d\phi^2(s)=0$.
\end{enumerate}
\end{definition}
It is well known that the above problem has a unique solution (\cite{cheman, Harrison-Reiman}). Define $\Gamma(\psi)=\phi^1$.
$\Gamma$ is referred to as {\it Skorohod map}.
In fact, $\Gamma$ has an explicit form given by
$$\displaystyle\Gamma(\psi)(t)= \psi(t)+\sup_{0\leq s\leq t} (\psi(s))^-. $$
It is easy to see that $\Gamma$ satisfies Lipschitz property i.e.,
\begin{equation}\label{skoro-lip}
 |\Gamma(\psi^1)-\Gamma(\psi^2)|^*_T\leq 2 |\psi^1-\psi^2|^*_T,
\end{equation}
for $\psi^i\in D([0,\iy), \R),\ \psi^i(0)\in\R_+,\ i=1,2.$

\section{Moderate deviations for many server queues}\label{sec2}
In this section, we introduce a single class G/M/N model. We consider a parallel
server system with single customer class and a pool of identical servers. We assume a buffer of infinite capacity. We parametrize the
system with parameter $n$. In what follows $n$ will be used to indicate the parametrization not as exponent, unless otherwise mentioned.
Let $\la^n$
 be given parameter where $\frac{1}{\la^n}$ represents the mean of the inter-arrival times of customers in the $n$-th system.
Let $\{IA(l) : l\in\N\}$ be a given sequence of i.i.d. of positive random variables with mean $\E[IA(1)]=1$ and variance
 $\text{Var}(IA(1))=\sigma^2_{IA}$ (here $IA$ stands for \textit{inter-arrival}). Assuming $\sum_{k=1}^0=0$, the number of arrivals of customers up to time $t$, for the
 $n$-th system, is given by
 $$ A^n(t)=\sup\Big\{ l\geq 0\ : \sum_{k=1}^l\frac{{\it IA}(k)}{\la^n}\leq t\Big\}, \quad t\geq 0.$$
Service time distributions are exponential. This justifies the use of notation M in G/M/N.
Let $\mu^n$ be the rate at which customers are served in the $n$-th system. 
N stands for the number of server which is assumed to vary with $n$. We denote by $N^n\in\N$
the number of servers in the $n-$th system.
We also consider the {\it moderate deviation rate} parameters $\{b_n\}$ with the property that $\lim b_n=\infty$ while
$\lim\frac{b_n}{\sqrt{n}}=0$.
 We assume that as $n\to\iy$,
\be\label{601}
\frac{\la^n}{n}\to\la\in (0 ,\iy),\quad \frac{N^n\mu^n}{n}\to\mu\in (0 ,\iy),\quad \frac{\sqrt{n}}{b_n}(\frac{\la^n}{n}-\frac{N^n\mu^n}{n})\to r\in(-\iy,\iy).
\ee
It is easy to see that under \eqref{601}, $\la=\mu$ i.e., the system is critically loaded. A similar condition in \cite{Puhal-1999}
is referred to as {\it near-heavy-traffic} condition.  We assume that each arriving customer has a single service requirement and it leaves the system when the job is completed by one of the servers otherwise it waits in the queue.

Let $X^n$ denote the number of customers in the system. Let $S(\cdot)$ be a standard Poisson process independent of the arrival process.
The number of service completion of jobs by time $t$ is given by
\be\label{603}
D^n(t)=S(\mu^n\int_0^tZ^n(s)),
\ee
where
$Z^n(t)$ denote the number of customers in service(or being served) at time $t$.
Hence we have
\be\label{605}
X^n(t)=X^n(0)+A^n(t)-D^n(t).
\ee
The system is assumed to work under non-idling policy i.e., $Z^n=X^n\wedge N^n$.
Next we define the scaled processes as follows
\be\label{607}
 \tilde{A}^n(t)=\frac{1}{b_n\sqrt{n}}(A^n(t)-\la^n t),
  \quad \tilde{S}_{\mu}^n(t)=\frac{1}{b_n\sqrt{n}}(S(N^n\mu^nt)-N^n\mu^n t),
  \ \tilde{X}^n(t)=\frac{1}{b_n\sqrt{n}}(X^n(t)-N^n).
\ee
It is easy to see from \eqref{605} that
\begin{align}
 \tilde{X}^n(t) &=\tilde{X}^n(0) + y^nt
 +\tilde{A}^n(t)-\tilde{S}_{\mu}^n(\frac{1}{N^n}\calT^n(t))
 + \frac{N^n\mu^n}{n}\frac{\sqrt{n}}{b_n}(t-\frac{1}{N^n}\calT^n(t))\notag
 \\
 &=\tilde{X}^n(0) + y^nt
 +\tilde{A}^n(t)-\tilde{S}_{\mu}^n(\frac{1}{N^n}\calT^n(t))
 + \frac{N^n\mu^n}{n}\frac{n}{N^n}\int_0^t(\tilde X^n(s))^-ds,\label{608}
 \end{align}
where $y^n=\frac{\sqrt{n}}{b_n}(\frac{\la^n}{n}-\frac{N^n\mu^n}{n}),\ \calT^n(t)=\int_0^tZ^n(s)ds$. We fix $T>0$ and assume:

\begin{condition}\label{or-mod}
The process $(\tilde A^n, \tilde S^n_\mu)$ satisfies large deviation principle (LDP) in $D([0, T],\R^2)$ with parameter $b_n^2$ and rate function $\Ir$ that takes value $\infty$ on discontinuous paths.
\end{condition}

\begin{remark}
Because of independence, it is enough if the processes $A^n, S^n_\mu$ satisfy LDP individually. In fact, one can impose some sufficient conditions
on the inter-arrival processes so that Condition \ref{or-mod} holds (see Remark \ref{rem1} and \ref{rem2} below).
\end{remark}
We also assume that the initial condition is deterministic and
$$\tilde X^n(0)\to x\in\R,\quad \text{as}\quad n\to\iy.$$
We are interested to find the rate function for $\tilde X^n$. We subdivide the problem in two theorems.

\spa

\begin{theorem}\label{main3}
Assume Condition \ref{or-mod} holds and $N^n=n$. Then $\{\tilde X^n\}_{n\geq 1}$ defined in \eqref{607} satisfies LDP in $D([0,T], \R)$ with parameter $b^2_n$ and rate
function $\Ir_X$ given by
$$ \Ir_X(\psi)=\inf_{\psi=G(\tilde\psi^1,\tilde\psi^2)}\Ir(\tilde\psi^1,\tilde\psi^2), $$
where $G(\tilde\psi^1,\tilde\psi^2)$ denotes the solution to the equation
\be
\psi= x + r t
 +\tilde\psi^1(t)-\tilde\psi^2(t)
 + \mu\int_0^t(\psi(s))^-ds.\label{614}
\ee
\end{theorem}

\noi{\bf Proof:}\ From \eqref{608}, we have
\be\label{609}
\tilde X^n(t)=\tilde{X}^n(0) + y^nt
 +\tilde{A}^n(t)-\tilde{S}_{\mu}^n(\frac{1}{N^n}\calT^n(t))
 + \frac{N^n\mu^n}{n}\int_0^t(\tilde X^n(s))^-ds.
\ee
Now given any tuple $(\tilde x, \tilde y, \kappa, \tilde\psi^1,\tilde\psi^2)\in\R^3\times D([0,T], \R^2)$ it is
easy to see that there exists a unique $\xi\in D([0, T], \R)$ satisfying the following:
\begin{align}
\xi(t)&=\tilde x + \tilde y t
 +\tilde\psi^1(t)-\tilde\psi^2(t)
 + \kappa\int_0^t(\xi(s))^-ds,\label{613}
 \\
|\xi|^*_T & \leq e^{\kappa T}\Big(|\tilde x+\tilde y t+\tilde\psi^1(t)-\tilde\psi^2(t)|^*_T \Big).\label{610}
\end{align}
Since $(\tilde A^n, \tilde S^n_\mu)$ satisfies LDP with rate function $\Ir$ and $\{\Ir\leq a\}, a\geq 0,$ is compact, we have
\be\label{611}
\lim_{\al\to\iy}\limsup_{n\to\iy}\p^{\frac{1}{b^2_n}}(|\tilde A^n|^*_T+|\tilde S^n_\mu|^*_T|\geq \al)=0.
\ee
Now for any $\del>0$,
$$\p(|t-\frac{1}{N^n}\calT^n(t)|^*_T>\del)\leq \p(e^{\frac{N^n\mu^n}{n}T}|\tilde X^n(0)+y^nt+\tilde A^n-\tilde{S}_{\mu}^n(\frac{1}{N^n}\calT^n(t))|^*_T\geq\frac{\sqrt n}{b_n}\frac{\del}{T}),$$
where we have used \eqref{610}. Therefore using \eqref{601} and \eqref{611}, we have
\be\label{612}
\limsup_{n\to\iy}\p^{\frac{1}{b^2_n}}(|t-\frac{1}{N^n}\calT^n(t)|^*_T>\del)=0.
\ee
Hence the sequence $\{\frac{1}{N^n}\calT^n\}$ converges super-exponentially in probability at rate $\frac{1}{b_n^2}$ to $\e(\cdot)$ in $D([0,T], \R)$ where $\e(t)=t$.
Therefore $(\tilde A^n, \tilde S^n_\mu\circ (\frac{1}{N^n}\calT^n))$ satisfies LDP with rate function $\Ir$ in $D([0,T],\R^2)$ (\cite{Puhal-Whit}, Lemma 4.3).
Denote $\xi$ by $G(\tilde x, \tilde y, \tilde \psi^1, \tilde\psi^2, \kappa)$ where $\xi$ satisfies \eqref{613}. Let $(\tilde\psi^1_n,\tilde\psi^2_n)\to (\tilde\psi^1,\tilde\psi^2)$
for some continuous path $(\tilde\psi^1,\tilde\psi^2)$. Let $\xi^n=G(\tilde X^n(0), y^n, \tilde \psi^1_n, \tilde\psi^2_n, \kappa^n)$, $\kappa^n= \frac{N^n\mu^n}{n}$, and
 $\xi=G(x, r, \tilde \psi^1, \tilde\psi^2, \mu)$. Then it is easy to see that $|\xi^n-\xi|^*_T\to 0$ as $n\to \infty$. Therefore extended contraction mapping principle yields that
 $X^n$ satisfies LDP with parameter $b_n^2$ and rate function
 $$ \Ir_X(\psi)=\inf_{\psi=G(\tilde\psi^1,\tilde\psi^2)}\Ir(\tilde\psi^1,\tilde\psi^2),$$
 where $G(\tilde\psi^1,\tilde\psi^2)$ denotes the solution to \eqref{614}.\hfill $\Box$

 \begin{theorem}\label{main4}
Assume Condition \ref{or-mod} holds. Let $N^n=o(n)$ and $x\in\R_+$. Then $\{\tilde X^n\}_{n\geq 1}$ defined in \eqref{607} satisfies LDP in $D([0,T], \R)$ with parameter $b^2_n$ and rate
function $\bar\Ir_X$ given by
$$ \bar\Ir_X(\psi)=\inf_{\psi=\Gamma(x+r\e+\tilde\psi^1-\tilde\psi^2)}\Ir(\tilde\psi^1,\tilde\psi^2), $$
where $\Gamma(\cdot)$ denotes the Skorohod map.
 \end{theorem}

\noi{\bf Proof:}\ From \eqref{608}, we have
\be\label{615}
\tilde X^n(t)=\tilde{X}^n(0) + y^nt
 +\tilde{A}^n(t)-\tilde{S}_{\mu}^n(\frac{1}{N^n}\calT^n(t))
 + \frac{N^n\mu^n}{n}\frac{n}{N^n}\int_0^t(\tilde X^n(s))^-ds.
\ee
By our assumption on $N^n$, we have $\frac{n}{N^n}\to\infty$ as $n\to\iy$. Given $\del>0$, we define the $\del$-oscillation function $osc_\del: D([0, T], \R)\to \R_+$
as follows:
$$osc_\del(\psi)=\sup\{|\psi(t)-\psi(s)|\ :\ |t-s|\leq\del\}.$$
By Condition \ref{or-mod} for any $\alpha>0$, $\{\Ir\leq\alpha\}$ is a compact set of continuous paths on $[0, T]$ to $\R^2$.
Therefore, given any $\del>0$, we have $\del_1>0$ so that  $osc_{\del_1}(\psi^1)+osc_{\del_1}(\psi^2)<\del$ for all $(\psi^1, \psi^2)\in \{\Ir\leq\alpha\}$. Therefore
\be\label{616}
\limsup_{\del_1\to 0}\limsup_{n\to\iy}\p^{\frac{1}{b_n^2}}(osc_{\del_1}(\tilde A^n)+osc_{\del_1}(\tilde S^n_\mu)\geq\del)=0.
\ee
Now choose $\eps>0$. We claim that
\be\label{617}
\limsup_{n\to\iy}\p^{\frac{1}{b_n^2}}(|(\tilde X^n)^-|^*_T\geq 3\eps)=0.
\ee
Define $\Om^n:=\{|(\tilde X^n)^-|^*_T\geq 3\eps\}$. Choose $n$ large enough so that $\tilde X^n(0)\geq -\frac{\eps}{2}$. For each $\om\in\Om^n$,
we will have random times $0\leq\sigma^n_1<\sigma^n_2\leq T$ such that $\tilde X^n(\sigma^n_1)> -\frac{3\eps}{2}, \tilde X^n(\sigma^n_2)\leq -\frac{5\eps}{2},$ and
$\tilde X^n(s)\leq -\eps$ on $[\sig^n_1, \sig^n_2]$. This is possible to do as the jump size of $\tilde X^n$ is $\frac{1}{b_n\sqrt n}$. Hence from \eqref{615}, we have
\be\label{618}
-\eps\geq y^n(\sig^n_2-\sig^n_1)-osc_{\sig^n_2-\sig^n_1}(\tilde A^n)-\osc_{\sig^n_2-\sig^n_1}(\tilde S^n_\mu)+\frac{N^n\mu^n}{n}\frac{n}{N^n}\eps(\sig^n_2-\sig^n_1).
\ee
Now if $(\sig^n_2-\sig^n_1)\geq\del_1$ for some fix $\del_1>0$ then \eqref{618} implies that
$2(|y^n|T+|\tilde A^n|^*_T+|\tilde S^n_\mu|^*_T)\geq \frac{N^n\mu^n}{n}\frac{n}{N^n}\eps\del_1$. If $(\sig^n_2-\sig^n_1)<\del_1$ then $|y^n\del_1|+osc_{\del_1}(\tilde A^n)+osc_{\del_1}(\tilde S^n_\mu)\geq \eps$. Therefore if we choose $\del_1>0$ so that $|y^n\del_1|<\frac{\eps}{2}$ for all $n$ large, then
$$\p (|(\tilde X^n)^-|^*_T\geq 3\eps)\leq \p(|\tilde A^n|^*_T+|\tilde S^n_\mu|^*_T\geq \kappa(\frac{n}{N^n})\eps\del_1)
+\p(osc_{\del_1}(\tilde A^n)+osc_{\del_1}(\tilde S^n_\mu)\geq\eps/2),$$
where $\kappa(\frac{n}{N^n})\to\iy$ as $n\to\iy$. Therefore first letting $n\to\iy$ and then letting $\del_1\to 0$ and
using \eqref{611} and \eqref{616}, the claim \eqref{617} follows.
Now rewriting \eqref{615} as
\be\label{619}
(\tilde X^n(t))^+=(\tilde X^n(t))^-+\tilde{X}^n(0) + y^nt
 +\tilde{A}^n(t)-\tilde{S}_{\mu}^n(\frac{1}{N^n}\calT^n(t))
 + \frac{N^n\mu^n}{n}\frac{n}{N^n}\int_0^t(\tilde X^n(s))^-ds,
\ee
we see that $((\tilde X^n(\cdot))^+, \frac{N^n\mu^n}{n}\frac{n}{N^n}\int_0^\cdot(\tilde X^n(s))^-ds)$ solves Skorohod problem for the date $(\tilde X^n(t))^-+\tilde{X}^n(0) + y^nt
 +\tilde{A}^n(t)-\tilde{S}_{\mu}^n(\frac{1}{N^n}\calT^n(t))$. Hence using the Lipschitz property of the Skorohod map \eqref{skoro-lip}
 we have
\be\label{620}
|\frac{N^n\mu^n}{n}\frac{n}{N^n}\int_0^\cdot(\tilde X^n(s))^-ds|^*_T
\leq 2|(\tilde X^n(t))^-+\tilde{X}^n(0) + y^nt
 +\tilde{A}^n(t)-\tilde{S}_{\mu}^n(\frac{1}{N^n}\calT^n(t))|^*_T.
\ee
 Since $\frac{\sqrt n}{b_n}(t-\frac{1}{N^n}\calT^n(t))=\frac{n}{N^n}\int_0^t(\tilde X^n(s))^-ds$, applying \eqref{611}, \eqref{617} and \eqref{620}, we have
 $$\limsup_{n\to\iy}\p^{\frac{1}{b^2_n}}(|t-\frac{1}{N^n}\calT^n(t)|^*_T>\del)=0. $$
Hence
$(\tilde A^n, \tilde S^n_\mu\circ (\frac{1}{N^n}\calT^n))$ satisfies LDP with rate function $\Ir$ (\cite{Puhal-Whit}, Lemma 4.3).
Now we consider a sequence $(\tilde\psi^1_n,\tilde\psi^2_n)\to (\tilde\psi^1,\tilde\psi^2)$ as $n\to\iy$
for some continuous path $(\tilde\psi^1,\tilde\psi^2)\in D([0, T],\R^2)$. Let $\xi^n$ be the solution to \eqref{613} with the data
$(\tilde X^n(0), y^n, \frac{N^n\mu^n}{n}\frac{n}{N^n}, \tilde\psi^1_n,\tilde\psi^2_n)$. Let $\xi=\Gamma(x+y\e+\tilde\psi^1-\psi^2)$.
To complete the proof it is enough to show that $|\xi^n-\xi|^*_T\to 0$ as $n\to\iy$. The proof will follow from the extended contraction mapping principle.
Given $\eps>0$, we choose $\del>0$ such that $(\om_\del(\psi^1_n)+\om_\del(\psi^2_n))<\frac{\eps}{4}$ for all $n$ large.
Since $\sup_n(|\tilde\psi^1_n|^*_T+|\tilde\psi^2_n|^*_T)<\iy$, we can choose $\del$ small enough to conclude that
$$|(\xi^n)^-|^*_T\leq 3\eps,$$
for large $n$ (using \eqref{618}). From \eqref{613}, we note that $(\xi^n)^+$ solves Skorohod problem for the data
$(\xi^n(\cdot))^-+\tilde{X}^n(0) + y^n\e(\cdot) +\tilde\psi^1_n(\cdot)-\tilde\psi^2_n(\cdot)$ and therefore Lthe ipschitz property
of the Skorohod map \eqref{skoro-lip} implies
\begin{eqnarray*}
|(\xi^n)^+-\xi|^*_T
&\leq & 2|(\xi^n(\cdot))^-+\tilde{X}^n(0) + y^n\e(\cdot) +\tilde\psi^1_n(\cdot)-\tilde\psi^2_n(\cdot)-
(x + r\e(\cdot) +\tilde\psi^1(\cdot)-\tilde\psi^2(\cdot))|^*_T,
\\
&\leq & 8\eps,
\end{eqnarray*}
for all $n$ large. Hence $|\xi^n-\xi|^*_T\to 0$ as $n\to\iy$. This completes the proof.\hfill $\Box$

\section{Control of multi-class G/M/N}\label{sec3}
 In this section, we introduce a multi-class G/M/N model and a related control problem. We consider a parallel
 server system with $\bI$ number of customer classes and a pool of identical servers. Let $\I=\{1,2,\ldots, \bI\}$. Let $\la_i^n>0, n\in\N, i\in\I,$
 be given parameter where $\frac{1}{\la_i^n}$ represents the mean of the inter-arrival time of class-$i$ customers in the $n$-th system.
 Given are $\bI$ independent sequence of i.i.d. $\{IA_i(l) : l\in\N\}_{i\in\I}$ of positive random variables with mean $\E[IA_i(1)]=1$ and variance
 $\text{Var}(IA_i(1))=\sigma^2_{i,IA}$. Assuming $\sum_{k=1}^0=0$, the number of arrivals of class-$i$ customers up to time $t$, in the
 $n$-th system, is given by
 $$ A_i^n(t)=\sup\Big\{ l\geq 0\ : \sum_{k=1}^l\frac{{\it IA}_i(k)}{\la_i^n}\leq t\Big\}, \quad t\geq 0.$$
$N^n\in\N$ denotes the number of servers in the $n-$th system. Service time distributions are exponential, with class dependent parameter.
Let $\mu^n_i$ be the rate at which class$-i$ customers are served in the $n$-th system.
We also consider the {\it moderate deviation rate} parameters $\{b_n\}$ with the property that $\lim b_n=\infty$ while
$\lim_{n\to\infty}[\frac{b_n}{\sqrt{n}}\vee\frac{N^n}{n}]=0$. Note that $N^n=o(n)$.
 We assume that as $n\to\iy$,
\begin{itemize}
\item $\frac{\la_i^n}{n}\to\la_i\in (0, \iy)$ and $\frac{N^n\mu_i^n}{n}\to\mu_i\in (0, \iy)$,
\item $\tilde{\la}_i^n:=\frac{1}{b_n\sqrt{n}}(\la_i^n-n\la_i)\to\tilde\la_i\in(-\infty, \infty)$,
\item $\tilde\mu_i^n:=\frac{1}{b_n\sqrt{n}}(N^n\mu_i^n-n\mu_i)\to\tilde\mu_i\in(-\infty, \infty)$.
\end{itemize}
Hence the traffic intensity for class-$i$, namely $\frac{\la_i^n}{N^n\mu_i^n}$, has limit $\rho_i:=\frac{\la_i}{\mu_i}$. The system
is assumed to be critically loaded i.e., $\sum_{i=1}^{\bI}\rho_i=1$.

\spa

\noi Let $B^n_i(t)$ be the number of servers working on class-$i$ customers at time $t\geq 0$. Therefore $B^n=(B^n_1,\ldots,B^n_{\bI})$ takes value in $(\{0\}\times\N)^\bI$.
Let $X^n_i, Q^n_i, I^n$ denote the number of class-$i$ customers in the system, the queue length of class-$i$ customers in the buffer and the number of
servers that are idle, respectively. Hence we have

\begin{align}
\displaystyle X_i^n &= Q^n_i+ B^n_i, \quad i\in\I,\label{1}
\\
\displaystyle N^n &= I^n+\sum_{i\in\I}B^n_i.\label{2}
 \end{align}
We are given $\bI$ independent standard Poisson processes $S_i, i\in\I$. The number of service completion of class$-i$ jobs by time $t$ is given by
\be\label{3}
D^n_i(t)=S_i(\mu^n_i\calT^n_i(t)),
\ee
where
\be\label{4}
\calT^n_i=\int_0^tB^n_i(s)ds.
\ee
Hence we have
\be\label{5}
X^n_i(t)=X_i^n(0)+A^n_i(t)-D^n_i(t).
\ee
For simplicity, the initial condition $X^n(0)=(X^n_i(0),\ldots,X^n_{\bI}(0))$ is assumed to be deterministic. The processes  $A^n, X^n, Q^n, B^n$ will always be assumed to have RCLL sample paths.
We will also assume that the processes $A^n_i, S_i, i\in\I,$ are mutually independent.

The process $B^n$ is regarded as control, that is determined based on the observation from the past (and present) events in the system. Fix $T>0$. Given $n$, the process $B^n$ is said
to be an {\it admissible control} if its sample paths lie in $D([0, T], \R^{\bI}_+)$ and
\begin{itemize}
\item $B^n(t)\in(\{0\}\times \N)^{\bI}$ for all $t\geq 0$;
\item For $i\in\I$ and $t\geq 0$,
\be\label{6}
B^n_i(t)\leq X^n_i, \quad \text{and}\quad \sum_{i\in\I}B^n_i(t)\leq N^n;
\ee
\item It is adapted to the filtration
$$\sigma\{A^n_i(s), D^n_i(s), i\in\I, s\leq t\}.$$
\end{itemize}
Denote the class of all admissible controls $B^n$ by $\frB^n$. We can see that under admissible control each server is allowed to serve a single customer at a time. We do not allow processor sharing. It is also easy to see that $\frB^n$ is non empty. For instance, if we define $B^n=0$ then $B^n\in\frB^n$.


Next we introduce the scaled processes. For $i\in\I$, let
\begin{align}\label{7}
 \tilde{A}^n_i(t) &=\frac{1}{b_n\sqrt{n}}(A_i^n(t)-\la^n_i t),
  \quad \tilde{S}_{\mu_i}^n(t)=\frac{1}{b_n\sqrt{n}}(S_i^n(N^n\mu^n_it)-N^n\mu^n_i t),\notag\\
  \ \tilde{X}_i^n(t)&=\frac{1}{b_n\sqrt{n}}(X_i^n(t)-\rho_iN^n).
\end{align}
It is easy to check from \eqref{5} that

\be\label{8}
 \tilde{X}_i^n(t)=\tilde{X}_i^n(0) + y_i^nt
 +\tilde{A}_i^n(t)-\tilde{S}_{\mu_i}^n(\frac{1}{N^n}\calT_i^n(t))
 + Z^n_i(t),
 \ee
where we denote
\be\label{9}
Z^n_i(t)=\frac{N^n\mu_i^n}{n}\frac{\sqrt{n}}{b_n}(\rho^i t-\frac{1}{N^n}\calT_i^n(t)),
\qquad y^n_i=\tilde\la_i^n-\rho_i\tilde\mu_i^n.
\ee
Since $\sum_{i\in\I}B^n_i\le N^n$ and $\sum_i\rho_i=1$, we see that
\be\label{10}
\sum_i\frac{n}{N^n\mu^n_i}Z^n_i \quad \text{starts from zero and is nondecreasing,}
\ee
The initial condition $X^n_i(0)$ is assumed to satisfy the following:
$$ \tilde{X}^n_i(0)\to x\in \R^{\bI}_+, \quad \text{as}\ n\to\infty.$$

The scaled arrival processes $\tilde A^n$ is assumed to satisfy a {\it moderate deviation
principle}. Let us first define the rate functions. Let
$\Ir_k, k=1,2,$ be functions on $D([0,T],\R^{\bI})$ defined as follows. For
$\psi=(\psi_1,\ldots,\psi_{\bI})\in D([0,T],\R^{\bI})$,
\[\Ir_1(\psi)=\left\{\begin{array}{ll}
               \frac{1}{2}\sum_{i=1}^\bI\frac{1}{\la_i\sig^2_{i,IA}}
               \int_0^T\dot\psi_i^2(s)ds &\ \mbox{if all}\ \psi_i\ \mbox{are absolutely continuous and}\ \psi(0)=0,
\\
                \infty & \ \mbox{otherwise},
              \end{array}
\right.\]
and
\[\Ir_2(\psi)=\left\{\begin{array}{ll}
               \frac{1}{2}\sum_{i=1}^\bI\frac{1}{\mu_i}
               \int_0^T\dot\psi_i^2(s)ds &\ \mbox{if all}\ \psi_i\ \mbox{are absolutely continuous and}\ \psi(0)=0,
\\
                \infty & \ \mbox{otherwise}.
              \end{array}
\right.\]

\begin{condition}\label{moderate} {\bf (Moderate deviation principle)}
The sequence
$\tilde{A}^n=(\tilde{A}_1^n,\ldots,
\tilde{A}_I^n),$
satisfies the LDP with parameters $b^2_n$ and rate function
$\Ir_1$ in $D([0,T],\R^{\bI})$; i.e.,
\begin{itemize}
 \item For any closed set $F\subset D([0,T],\R^{\bI})$
$$\limsup\frac{1}{b_n^2}\log\p(\tilde{A}^n\in F)\leq -\inf_{\psi\in F}\Ir_1(\psi),$$
\item For any open set $G\subset D([0,T],\R^{\bI})$
$$\liminf\frac{1}{b_n^2}\log\p(\tilde{A}^n\in G)\geq -\inf_{\psi\in G}\Ir_1(\psi).$$
\end{itemize}
\end{condition}
\begin{remark}\label{rem1}
It is shown in \cite{Puhal-Whit} that each one of the following statements is sufficient for
Condition \ref{moderate} to hold:
\begin{itemize}
\item There exist constants $a_0>0$, $\beta\in(0,1]$ such that $E[e^{a_0(IA_i)^\beta}]<\iy$ , $i\in\I$, and $b_n^{\beta-2}n^{\beta/2}\to\infty$;
\item For some $\del>0$, $E[(IA_i)^{2+\del}]<\infty$, $i\in\I$, and
    $b_n^{-2}\log n\to\infty$.
\end{itemize}
\end{remark}

\begin{remark}\label{rem2}
Since the inter-arrival time for a Poisson process is exponential,
using Remark \ref{rem1}, we see that $\tilde S^n_\mu=(\tilde S^n_{\mu_1},\ldots,\tilde S^n_{\mu_{\bI}})$,
satisfies Large deviation principle in $D([0,T],\R^{\bI})$ with parameter $b_n^2$ and rate function $\Ir_2$.
Therefore using the independence of the processes (see \cite{Lynch-Sethu}) and extended contraction mapping principle we see that
$(\tilde A^n, \tilde S^n_\mu)$ satisfies Large deviation principle in $D([0,T],\R^{2\bI})$ with parameter $b^2_n$
and rate function $\Ir(\psi)=\Ir_1(\psi^1)+\Ir_2(\psi^2)$, $\psi=(\psi^1,\psi^2)\in D([0,T],\R^{2\bI})$.
\end{remark}

To present our control problem, we consider nonnegative functions $h$ and $g$ from $\R^{\bI}$ to $\R$ which are nondecreasing with respect to the usual
partial order on $\R^{\bI}$. We assume that $h, g$ have at most linear growth, i.e., there exist constants $C_1, C_2$ such that
$$ g(x)+h(x)\leq C_1\|x\|+C_2.$$
Given $n$, the cost associated with the initial condition $\tilde X^n(0)$ and control $B^n$ is given by
\[
J^n_X(\tilde X^n(0), B^n)=\frac{1}{b_n^2}
\log\E\Big[e^{b^2_n[\int_0^Th(\tilde{X}^n(s))ds+g(\tilde{X}^n(T))]}\Big].
\]
We are interested to analyze the value function
$$V^n_X(\tilde X^n(0))=\inf_{B^n\in\frB^n} J^n_X(\tilde{X}^n(0),B^n).$$
We now introduce another value function associated to the queue length. To do this, we define $\tilde Q^n_i=\frac{1}{b_n\sqrt{n}}Q^n_i, i\in\I,$ and $\tilde Q^n=(\tilde Q^n_1,\ldots,\tilde Q^n_{\bI})$. Let
\[
J^n_Q(\tilde Q^n(0), B^n)=\frac{1}{b_n^2}
\log\E\Big[e^{b^2_n[\int_0^Th(\tilde{Q}^n(s))ds+g(\tilde{Q}^n(T))]}\Big].
\]
The associated value function is given by
$$V^n_Q(\tilde Q^n(0))=\inf_{B^n\in\frB^n} J^n_Q(\tilde{X}^n(0),B^n).$$

\subsection{A differential game and main results}\label{sec3.1}

We next develop a differential game for the limiting behavior
of the value functions defined above. This game problem has been studied in \cite{Atar-Bis}.
Let $\theta=(\frac{1}{\mu_1},\ldots,\frac{1}{\mu_{\bI}})$ and $y=(y_1,\ldots,y_{\bI})$
where $y_i=\tilde\la_i-\rho_i\tilde\mu_i$.
Denote $P=C_0([0,T],\Ri)$ (the subset of $C([0,T],\Ri)$ of functions with initial value $0$) and
\[
E=\{\zeta\in C([0,T],\R^{\bI}):\theta\cdot\zeta \text{ starts from zero and is nondecreasing}\}.
\]
The topology on both the spaces are induced by uniform topology.
Let $R$ be a mapping from $D([0,T],\R^{\bI})$ into itself defined by
\begin{equation}\label{2001}
R_i[\psi](t)=\psi_i(\rho_it),\qquad t\in[0,T],\, i\in\calI.
\end{equation}
Given $\psi=(\psi^1,\psi^2)\in P$ and $\zeta\in E$, we define the {\it dynamics associated
with initial condition $x$ and data $\psi,\zeta$} as
\be
\varphi_i(t)= x_i+y_it+ \psi^1_i(t)-R_i[\psi^2](t)+\zeta_i(t),\qquad i\in\calI.
\label{201}
\ee
It is easy to see the analogy between the above equation and equation \eqref{8},
and between the condition $\theta\cdot\zeta$ nondecreasing and property \eqref{10}.
The following condition will also be used,
\be\label{13}
\varphi_i(t)\ge0,\qquad t\ge0,\ i\in\calI.
\ee
To define the game in the sense of Elliott and Kalton \cite{Elli-Kal},
we need the notion of strategies.
A measurable mapping $\alpha:P\to E$ is called a {\it strategy for the minimizing player}
if it satisfies the causality property. Namely,
for every $\psi=(\psi^1,\psi^2), \tilde \psi=(\tilde\psi^1,\tilde\psi^2) \in P$
and $t\in[0,T]$,
\be\label{14}
\text{$(\psi^1,R[\psi^2])(s)=(\tilde\psi^1,R[\tilde\psi^2])(s)$
for all $s\in[0,t]$\quad implies\quad $\alpha[\psi](s)=\alpha[\tilde\psi](s)$
for all $s\in[0,t]$.}
\ee
Given an initial condition $x$,
a strategy $\al$ is said to be {\it admissible} if, for $\psi\in P$ and
$\zeta=\al[\psi]$, the corresponding dynamics \eqref{201} satisfies the nonnegativity
constraint \eqref{13}. The set of all admissible strategies for the minimizing
player is denoted by $A_x$.
Given $x$ and $(\psi,\zeta)\in P\times E$, we define the cost by
$$
c(\psi,\zeta)=\int_0^Th(\varphi(t))dt+g(\varphi(T))-\Ir(\psi),
$$
where $\varphi$ is the corresponding dynamics given by \eqref{201} and $\Ir$ is given in Remark \ref{rem2}.
The value of the game is defined by
$$
V(x)=\inf_{\alpha\in A_x}\sup_{\psi\in P}c(\psi, \alpha[\psi]).
$$
One can also obtain a simpler, equivalent formulation of the above game (see Remark 2.2 in \cite{Atar-Bis}).
\subsubsection{Main results}
Before we state our main results, let us introduce two conditions that will be used to prove the results. For $\om\in\R_+$, define
$$
h^*(w)=\inf\{h(x) : x\in\R^{\bI}_+, \theta\cdot x=w\},
\quad
g^*(w)=\inf\{g(x) : x\in\R^{\bI}_+, \theta\cdot x=w\}.
$$
We impose the following condition.
\begin{condition}\label{mini}
{\bf (Existence of a continuous minimizing curve)}\ There exists a continuous map $f : \R_+\to\R^{\bI}_+$ such that for all $w\in\R_+$,
$$ \theta\cdot f(w)=w,
\quad
h^*(w)=h(f(w)),
\quad
g^*(w)=g(f(w)).
$$
\end{condition}
We refer  to \cite{Atar-Bis} for the examples of $h$ and $g$ satisfying above condition.
Similar condition is also used in \cite{atasol}, \cite{atagur},
where an analogous many-server model is treated in a diffusion regime. We comment in Remark \ref{end-rem} about weakening this assumption.

\begin{condition}\label{unbounded}{\bf (Exponential moments)}
Denote
$
\La_T(\psi^1)=\sum_{i=1}^{\bI}\sup_{[0,T]}|\psi^1_i(t)| .
$
Then for any constant $K$,
$$
\limsup_{n\to\iy}\frac{1}{b^2_n}\log\E[e^{b^2_nK\La_T(\tilde{A}^n)}]
<\infty.$$
\end{condition}
In view of Proposition 2.1 in \cite{Atar-Bis}, if there exists $a_0>0$ such that $\sup_{i\in\I}\E[e^{a_0\,IA_i}]<\infty$ then Condition \ref{unbounded} holds.

\begin{remark}\label{rem3}
If Condition \ref{unbounded} holds, then it is easy to see that
for any constant $K$,
$$
\limsup_{n\to\iy}\frac{1}{b^2_n}\log\E[e^{b^2_nK(\La_T(\tilde{A}^n)+\La_T(\tilde{S}^n_\mu))}]
<\infty.$$
\end{remark}

Now we are ready to state our main results.
\begin{theorem}\label{main1}
Let Conditions \ref{moderate}, \ref{mini} and \ref{unbounded} hold and $\lim_{n\to\iy}\frac{N^n}{b_n\sqrt n}=0$. Then
$$\lim V^n_X(\tilde{X}_n(0))=V(x).$$
\end{theorem}

\begin{theorem}\label{main2}
Let Conditions \ref{moderate}, \ref{mini} and \ref{unbounded} hold and $\lim_{n\to\iy}\frac{N^n}{b_n\sqrt n}=0$. Then
$$\lim V^n_Q(\tilde{Q}^n(0))=V(x).$$
\end{theorem}

\begin{theorem}\label{main5}
Let Conditions \ref{moderate} and \ref{mini} hold. If $h$ and $g$ are bounded, then
$$\lim_{n\to\iy} V^n_{X}(\tilde{X}^n(0))=V(x).$$
\end{theorem}

\subsection{Proof of Theorem \ref{main1}, \ref{main2} and \ref{main5}}\label{sec3.2}
\subsubsection{Lower bound}\label{lower}
Before we go in further details, let us mention a solution to the above game problem that
was obtained in \cite{Atar-Bis}. Recall the one-dimensional {\it Skorohod map} $\Gam$ from $D([0, T], \R)$ into itself. Given $\psi=(\psi^1, \psi^2)\in P$,
define
\[
\hat\psi(t)=x+yt+\psi^1(t)-\psi^2(t),\qquad t\in[0,T].
\]
We define
\be\label{018}
\hat\bzeta[\psi](t)=f(\hat\varphi_\theta[\psi](t))-\hat\psi(t),\qquad t\in[0,T],
\ee
where $\hat\varphi_\theta[\psi]=\Gam[\theta\cdot\hat\psi]$. Let us define $\bzeta[\psi^1,\psi^2]=\hat\bzeta[\psi^1,R[\psi^2]]$
where $R[\cdot]$ is given by \eqref{2001}.
In \cite[Proposition 3.1]{Atar-Bis}, it is proved that $\bzeta$ is a minimizing strategy for the game i.e.,
\be
V(x)=\sup_{\psi\in P} c(\psi, \bzeta[\psi]).
\label{305}
\ee
In fact, $\bzeta$ satisfies the following minimality property: for any $\alpha\in A_x$ with $(\psi, \tilde\zeta), \tilde\zeta=\alpha[\psi],$ satisfying \eqref{201} with the associated dynamics $\tilde\varphi$  we have
\be\label{10000}
h(\tilde\varphi(t))\geq h(\varphi(t)), \quad g(\tilde\varphi(T))\geq g(\varphi(T)), \quad t\geq 0,
\ee
where $\varphi$ is the dynamics associated to $(\psi, \bzeta[\psi])$ satisfying \eqref{201}. It is obvious from the definition that an analogous
minimality property holds for $\hat\bzeta$.

For $\kap>0$, we define
\begin{equation}\label{19}
D(\kappa)=\{\psi=(\psi^1,\psi^2)\in D([0,T], \R^{2\bI})
:\ \|\psi^1\|^*_T+\|\psi^2\|^*_T\leq\kappa\ \mbox{and}\ \bar\psi(0)\in\R^{\bI}_+\},
\end{equation}
where
$$\bar\psi(t)=x+yt+\psi^1(t)-R[\psi^2](t),\qquad t\in[0,T].$$
It is shown in \cite[Proposition 3.2]{Atar-Bis} that there exists constant $\gam_1,\gam_2$ such that
\be
 \|\hat\bzeta[\psi]\|^*_t\leq \gam_1\big(\La_t(\psi^1)+\La_t(\psi^2)\big)+\gam_2,\quad t\in [0,T],
 \label{311}
\ee
for all $(\psi^1,\psi^2)\in D([0,T], \R^{2\bI})$. Given a map $\varphi:[0, T]\to\R^k$
and a constant $\eta>0$, we define the $\eta$-oscillation of $\varphi$ as
$$ \osc_\eta(\varphi)=\sup\{\|\varphi(s)-\varphi(t)\|\ :\ |s-t|\leq\eta,\  s,t\in[0, T]\}.$$
Then for any given $\kap, \eps>0$ there exists $\del, \eta$ such that the followings hold:
For any $\psi,\tilde\psi\in D(\kap)$
\be
 \|\hat\bzeta[\psi]
 -\hat\bzeta[\tilde\psi]\|^*_T\le \eps\ \mbox{ if }\ \|\psi^1-\tilde\psi^2\|^*
 +\|\psi^2-\tilde\psi^2\|^*\leq\delta,
 \label{315}
\ee
and
\begin{equation}
 \osc_\eta(\hat\bzeta[\psi])\le \eps\ \mbox{ provided }\ \osc_\eta(\psi)\le\delta.
 \label{316}
\end{equation}

\begin{theorem}\label{lower1}
Assume Conditions \ref{moderate} and \ref{mini} hold. Then $\liminf V^n_X(\tilde{X}_n(0))\geq V(x)$.
\end{theorem}

\noi{\bf Proof:}\
The proof of the theorem follows from \cite{Atar-Bis} except some
suitable modifications. We add here some details for clarity and convenience of the readers.
Fix $\tilde\psi=(\tilde\psi^1,\tilde\psi^2)\in P$.
Recall metric $\m(\cdot,\cdot)$  on $D([0,T], \R^{2\bI})$ which induces the $J_1$ topology.
Define, for $r>0$,
$$
\A_r=\{\psi\in D([0,T], \R^{2\bI}) \ :\ \m(\psi,\tilde\psi)<r\}.
$$
Since $\tilde\psi$ is continuous, for any $r_1\in (0, 1)$ there exists $r, \eta>0$ such that
\begin{equation}
\psi\in \A_r \quad\text{implies}\quad \|\psi-\tilde\psi\|^*<r_1,\quad \text{osc}_\eta(\psi)<r_1.
\label{411}
\end{equation}
This can be done as for any $f\in\Ups$ (see Notations),
\begin{eqnarray*}
\|\psi(t)-\tilde\psi(t)\| &\leq & \|\psi(t)-\tilde\psi(f(t))\|+\|\tilde\psi(f(t))-\tilde\psi(t)\|,
\\
|f(t)-t|^*_T &\leq & T(e^{\|f\|^\circ}-1),
\end{eqnarray*}
and $\tilde\psi$ is uniformly continuous on $[0, T]$.
Define $\theta^n=(\frac{n}{N^n\mu_1^n},\frac{n}{N^n\mu_2^n},\ldots,\frac{n}{N^n\mu_{\bI}^n})$. Then
$\theta^n\to\theta$ as $n\to\infty$.
Now, given $0<\eps<1$, choose a sequence of policies $\{B^n\}$ such that
$$
V^n_X(\tilde{X}^n(0)) + \eps> J_X(\tilde{X}^n(0), B^n)
\ \mbox{and}\ B^n\in\frB^n\ \mbox{for all}\ n.
$$
Recall
\be
 \tilde{X}_i^n(t) = \tilde{X}_i^n(0)+y_i^nt+
 \tilde{A}^n_i(t)-\tilde{S}^n_{\mu_i}(\frac{1}{N^n}\calT_i^n(t))
 +Z_i^n(t),
 \label{413}
\ee
where
\be\label{53}
Z_i^n(t)=\frac{N^n\mu_i^n}{n}\frac{\sqrt{n}}{b_n}(\rho_i t-\frac{1}{N^n}\calT_i^n(t)),\qquad
\calT_i^n(t) = \int_0^tB_i^n(s)ds.
\ee
We claim that for all $n$ large and $(\tilde A^n, \tilde S^n_\mu)\in\A_r,$
\be\label{666}
\sup_{i\in\I}|(\tilde X^n_i)^-|^*_T\leq (6+c)r_1,
\ee
where $c=\sup_n\|y^n\|$.
To prove the claim, let us take $i\in\I$ such that $|(\tilde X^n_i)^-|^*_T> (6+c) r_1$ for infinitely many $n$. For large $n$, we have $\tilde X^n_i(0)>-r_1$.
Hence we have times $\sig^n_1<\sig^n_2\leq T$ such that $\tilde X^n_i(\sig^n_1)\geq- 2r_1,\ \tilde X^n_i(\sig^n_2)\leq- (c+5)r_1$ and $\tilde X^n_i(s)\leq -r_1$
for all $s\in[\sig^n_1, \sig^n_2]$. Therefore $\frac{B^n_i(s)-\rho_iN^n}{b_n\sqrt n}\leq \frac{X^n_i(s)-\rho_iN^n}{b_n\sqrt n}\leq -r_1 $ for all $s\in[\sig^n_1, \sig^n_2]$.
Hence using \eqref{413} we have
\be\label{445}
-(c+3)r_1\geq y^n_i(\sig^n_2-\sig^n_1) - 2osc_{\sig^n_2-\sig^n_1}(\tilde A^n, \tilde S^n_\mu)+\frac{N^n\mu_i^n}{n}\frac{n}{N^n}r_1(\sig^n_2-\sig^n_1).
\ee
Using \eqref{411} and the fact $\frac{n}{N^n}\to\iy$, we see that \eqref{445} leads to a contradiction for large $n$ if $(\sig^n_2-\sig^n_1)\geq r_1\wedge\eta$.
Again if $(\sig^n_2-\sig^n_1)\leq r_1\wedge\eta$, then \eqref{445} is contradicting to \eqref{411} as the right most term in \eqref{445} is non-negative.
This proves the claim \eqref{666}.

Given $G>0$, define
$$
\tau_n=\inf\{t\ge0:\theta^n\cdot Z^n(t)>G\}\w T
\equiv\inf\Big\{t\geq 0\ :\ \frac{\sqrt{n}}{b_n}\Big(t-\frac{1}{N^n}\sum_{i=1}^\bI \calT_i^n(t)\Big)>G\Big\}
\wedge T.
$$
It is possible to choose $\kappa_1>0$ such that for $(\tilde A^n, \tilde S^n_\mu)\in\A_r$ and $t>\tau_n$ (see (4.6) in \cite{Atar-Bis}),
\be\label{738}
\theta^n\cdot\tilde X^n(t)\geq -\kappa_1+G.
\ee
We note that for any positive $\kap_1, \kap_2$,
\be\label{10001}
\inf\{h(x): \theta\cdot x\geq \kap_1, \ x_i\geq -\kap_2\}=\inf\{h(x): \theta\cdot x= \kap_1, \ x_i\geq -\kap_2\}.
\ee
To see this we consider $x$ with $\theta\cdot x> \kap_1, \ x_i\geq -\kap_2,$ and multiply all its positive coordinate 
by $\beta$. Call the new point as $x_{(\beta)}$. Therefore $x_{(1)}=x$ and $h(x_{(\beta)})$ decreases with $\beta$. Hence we can choose $\beta<1$ such that
$h(x_{(\beta)})\leq h(x)$ and $\theta\cdot x_{(\beta)}=\kap_1, \ x_{(\beta)i}\geq -\kap_2$. This proves \eqref{10001}.
Now we note that given $\kap>0$, for large $n$ (so that $\theta^n_i\leq2\theta_i, \ i\in \calI,$)
\be\label{001}
\theta^n\cdot x\geq\kap \To \theta\cdot x\geq \frac{\kap}{2}+\frac{\theta^n\cdot x^-}{2}-\theta\cdot x^-\geq \frac{\kap}{2}
-\theta\cdot x^-.
\ee
Next, let $\varphi:[0,T]\to\R^\bI$ be the dynamics corresponding to $(\tilde\psi,\zeta)$,
where $\zeta=\bzeta[\tilde\psi]$, namely
\be\label{12}
\varphi_i(t)=x_i+y_it+\tilde\psi^1_i(t)-R[\tilde\psi^2]_i(t)
+\zeta_i(t).
\ee
Then $\varphi(t)=f(\varphi_\theta[\tilde\psi](t))$ \eqref{018} where $\varphi_\theta[\tilde\psi]=\Gamma[\theta\cdot\bar\psi],
\bar\psi(t)=x+yt+\tilde\psi^{1}(t)-R[\tilde\psi^{2}](t)$.
Let $\bar\om_h\ [\bar\om_g]$ be the modulus of continuity of $h$ [resp. $g$] over
$\{x\in\R^{\bI}:\ x\cdot\theta\leq|\varphi_\theta[\tilde\psi]|^*_T, \ x_i\geq -(c+6)\}$.
For $\kappa,\ \frac{\kappa}{2}-\sqrt{I}\|\theta\|(c+6)>|\varphi_\theta[\tilde\psi]|^*_T,$
and large $n$ we have
\begin{eqnarray*}
\inf\{h(x) : \theta^n\cdot x\geq \kappa, x_i\geq -(c+6)r_1\}
&\geq & \inf\{h(x) : \theta\cdot x\geq |\varphi_\theta[\tilde\psi]|^*_T, x_i\geq -(c+6)r_1\}
\\
&\geq & \inf\{h(x) : x\in R^\bI_+,\ \theta\cdot x\geq |\varphi_\theta[\tilde\psi]|^*_T\}-\bar\om_h(\sqrt I(c+6)r_1)
\\
&\geq& |h(\varphi)|^*_T-\bar\om_h(\sqrt I(c+6)r_1),
\end{eqnarray*}
where for the first inequality we use \eqref{001}, for second inequality we use \eqref{10001}, and for the last inequality we use the monotonicity of $h$.
Similar estimate holds for $g$.
Therefore we can find $G>0$ such that
for all large $n$, $(\tilde A^n, \tilde S^n_\mu)\in\A_r$,
\begin{equation}
h(\tilde{X}^n(t))\geq |h(\varphi)|^*-\del(r_1)\ \mbox{and}\ \ g(\tilde{X}_n(t))
\geq g(\varphi(T))-\del(r_1),\label{444}
\end{equation}
on $\{t>\tau_n\}$ where $\del(r_1)\to 0$ as $r_1\to 0,$ and for $t\leq\tau_n$
\begin{equation}
 \|Z_n(t)\|\leq \kappa_2.
 \label{417}
\end{equation}
for some constant $\kappa_2$.
Hence using \eqref{413}, \eqref{411} and \eqref{417}, we obtain a constant $\kap_3$ such that for $(\tilde A^n, \tilde S^n_\mu)\in\A_r$
and all $n$ large
$$\sup_{t\in[0,\tau_n]}\|\tilde X^n(t)\|\leq\kap_3.$$
Consider the stochastic processes $Y^n, \tilde Y^n, \tilde Z^n$,
with values in $\R^{\bI}$ such that $(\tilde X^n_i)^+(t)=x_i+y_it+Y^n_i(t)-\tilde Y_i^n(t)+\tilde Z_i^n(t)$ on $[0,\tau_n]$
where
\begin{eqnarray*}
Y_i^n(t) &=& \tilde{A}^n_i(t\w\tau_n),
\\
\tilde{Y}_i^n(t) &=& x_i-\tilde{X}_i^n(0)-(\tilde X^n_i)^-(t)+(y_i-y_i^n)t+\tilde{S}^n_{\mu_i}(\frac{1}{N^n}\calT_i^n(t\w\tau_n))
-(1-\mu_i\theta^n_i)Z_i^n(t\w\tau_n),
\\
\tilde Z^n_i(t) &=& \mu_i\theta^n_i Z_i^n(t).
\end{eqnarray*}

Define $W^n(t)=x+yt+Y^n(t)-\tilde{Y}^n(t)+\hat\bzeta[Y^n,\tilde Y^n](t)$.

By \eqref{417}, we have $\sup_i\sup_{[0,\tau_n]}|\rho_it-\frac{1}{N^n}\calT^n_i(t)|\to 0$ as $n\to\infty$.
Therefore using the regularity property of $\hat\bzeta$ \eqref{316} and \eqref{411}, \eqref{666} with a proper choice of $r_1<\eps$ we have
\be
\\
\|\varphi-W_n\|^*_{\tau_n} \leq \kappa_4\eps.\label{419}
\ee
for all $n$ large and $(\tilde A^n, \tilde S^n_\mu)\in\A_r$ and some constant $\kappa_4$.
Denote by $\omega_h$ [$\omega_g$] the modulus of continuity of $h$ [resp., $g$]
over $\{q:\|q\|\leq \|\varphi\|^*_T+\kappa_3+\kap_4\}$. Hence using minimality property of $\hat\bzeta$ \eqref{10000} we have
\be\label{426}
\begin{cases}
h(\tilde{X}^n(t)) &\geq h((\tilde X^n_t)^+)-\om_h(\sqrt I(c+6)r_1)\geq h(W^n(t))-\om_h(\sqrt I(c+6)r_1),
\\
g(\tilde{X}^n(t))&\geq h(W^n(t))-\om_g(\sqrt I(c+6)r_1)\geq h(W^n(t))-\om_g(\sqrt I(c+6)r_1).
\end{cases}
\ee
for $t\in[0,\tau_n]$ and $(\tilde A^n, \tilde S^n_\mu)\in\A_r$ for large $n$.

Then by \eqref{419} and \eqref{426}, for $(\tilde A^n, \tilde S^n_\mu)\in\A_r$ and all large $n$,
\begin{eqnarray*}
 \int_0^{\tau_n}h(\tilde{X}_n(s))ds &\geq &\int_0^{\tau_n}h(W_n(s))ds-T\om_h(\sqrt I(c+6)\eps)
 \\
&\geq &\int_0^{\tau_n} h(\varphi(s))ds-T\omega_h(\kap_4\eps)-T\om_h(\sqrt I(c+6)\eps).
\end{eqnarray*}
Combined with \eqref{444} this gives
\[
\int_0^Th(\tilde X_n(s))ds\ge\int_0^Th(\varphi(s))ds-T\Big(\omega_h(\kap_4\eps)+\om_h(\sqrt I(c+6)\eps)+\del(\eps)\Big).
\]
A similar argument gives
\[
g(\tilde X_n(T))=
g(\varphi(T))\chi_{\{T\leq\tau_n\}}+g(\varphi(T))\chi_{\{T>\tau_n\}}
\ge g(\varphi(T))-\Big(\omega_g(\kap_4\eps)+\om_g(\sqrt I(c+6)\eps)+\del(\eps)\Big).
\]
Hence using condition \ref{moderate}, it is easy to show that for all $n$ large
\begin{eqnarray*}
\E[e^{b_n^2[\int_0^Th(\tilde{X}_n(s))ds+g(\tilde{X}_n(T))]}]
&\geq& \E\Big[e^{b_n^2[\int_0^Th(\tilde{X}_n(s))ds+g(\tilde{X}_n(T))]}
\chi_{\{(\tilde{A}_n, \tilde{S}_n)\in \A_r\}}\Big]
\\
&\geq & e^{b^2_n\Big( \int_0^Th(\varphi(s))ds+g(\varphi(T))-\Ir(\tilde\psi)-a(\eps)-\eps\Big)},
\end{eqnarray*}
and hence
\begin{eqnarray*}
\frac{1}{b_n^2}\log\E[e^{b_n^2[\int_0^Th(\tilde{X}_n(s))ds+g(\tilde{X}_n(T))]}]
\geq \int_0^Th(\varphi(s))ds+g(\varphi(T))-\Ir(\tilde\psi)-a(\eps)-\eps,
\end{eqnarray*}
where $a(\eps)=T\omega_h(\kap_4\eps)+T\om_h(\sqrt I(c+6)\eps)+\om_g(\sqrt I(c+6)\eps)+\omega_g(\kap_4\eps)+(T+1)\del(\eps)$.
The proof follows letting $n\to\infty$ first and then $\eps\to 0$.\hfill $\Box$

\begin{theorem}
Assume Conditions \ref{moderate} and \ref{mini} hold and $\lim\frac{N^n}{b_n\sqrt n}=0$. Then 
$$\liminf V^n_Q(\tilde{X}_n(0))\geq V(x).$$
\end{theorem}

\noi{\bf Proof:}\ We note that proof of  Theorem \ref{lower1} relies on two estimates, \eqref{444} and \eqref{419}. To get these estimates,
we note that for any $\kap>0$ and large $n$, $\theta^n\cdot\tilde X^n\geq \kap+1\To\theta\cdot\tilde Q^n\geq\kap/2$ and $\|\tilde X^n-\tilde Q^n\|^*_T=o(1)$.
Hence the proof follows.\hfill $\Box$

\subsubsection{Upper bound}

\begin{theorem}\label{upper1}
Assume Conditions \ref{moderate}, \ref{mini} and \ref{unbounded} hold and $\lim\frac{N^n}{b_n\sqrt n}=0$. Then
$$\limsup V^n_X(\tilde X^n(0))\leq V(x).$$
\end{theorem}

\begin{remark}
 If the functions $h, g$ are bounded then Condition \ref{unbounded} is not required
 in the above statement.
\end{remark}

The proof is based on the construction of a suitable admissible policy. The main idea of the proof is similar to that appear
in \cite{Atar-Bis}. However, the proof appear here is complicated than that appeared in \cite{Atar-Bis}. The main difficulty we face
here is due to the constrain on the policy that does not allow processor sharing.
 The idea is to make use of the preemptive behavior of the policy.
We construct a policy that serves each class of customers over small time intervals (defined in suitable sense) and on
average effort given to serve
 class $i$ is $\approx \rho_i+\mathfrak{E}$ where the correction $\mathfrak{E}$ is small and leads us to the correct limit.

\spa

\noi{\bf Proof:}\ Let $\Del>0$ be a given constant. Define
\be\label{60}
\calQ=\{\psi\in D([0, T], \R^{2\bI})\ : \ \Ir(\psi)\leq \Del\}.
\ee

By the definition of the rate function $\Ir$ (from Remark \ref{rem2}) , $\calQ$ is a compact set containing absolutely
continuous paths starting from zero (particularly, $\calQ\subset P$),
with derivative having $L^2$-norm uniformly bounded.
Consequently, there exists a constant $M=M_\Del$  such that
$\|\psi^1\|^*+\|\psi^2\|^*\leq M$ for all $\psi\in \calQ$. Consider the set $D(M+1)$ \eqref{19},
let $\eps\in(0,1)$ be given, and choose $\delta, \eta>0, \del\in (0, \eps),$
as in (\ref{315}) and (\ref{316}), corresponding to $\eps$ and $\kappa=M+1$.
It follows from the $L^2$ bound alluded to above, that for each fixed $\Del$,
the members of $\calQ$ are equicontinuous. Hence one can choose $v_0\in(0,\eta)$ (depending on
$\Del$), such that
\begin{equation}
 \osc_{v_0}(\psi^l_i)<\frac{\delta}{4\sqrt{2\bI}}, \ \mbox{for all}\ \psi=(\psi^1,\psi^2)\in \calQ,\,l=1,2,\,
 i\in\I.\label{456}
\end{equation}
Recall
$$\A_r(\tilde\psi)=\{\psi\in D([0,T], \R^{2I}) \ :\ \m(\psi,\tilde\psi)<r\}.$$
Noting that, for any $f\in\Ups$ (see Notations),
\begin{eqnarray*}
\|\psi(t)-\tilde\psi(t)\| &\leq & \|\psi(t)-\tilde\psi(f(t))\|+\|\tilde\psi(f(t))-\tilde\psi(t)\|,
\\
|f(t)-t|^*_T &\leq & T(e^{\|f\|^\circ}-1),
\end{eqnarray*}
it follows, by the equicontinuity of the members of $\calQ$, that it is possible to choose
$v_1>0$ such that, for any $\tilde\psi\in \calQ$,
\begin{equation}
\psi\in \A_{v_1}(\tilde\psi) \quad \text{ implies } \quad
\|\psi-\tilde\psi\|^*<\frac{\delta}{4}.\label{457}
\end{equation}
Let $v_2=\min\{v_0, v_1, \frac{\eps}{2}\}$. Since $\calQ$ is compact and $\Ir$ is lower semicontinuous,
one can find a finite number of members $\psi^1, \psi^2,\ldots,\psi^N$
of $\calQ$, and positive constants $v^1,\ldots, v^N$ with $v^k<v_2$, satisfying
$\calQ\subset\cup_k \A^k$, and
\begin{equation}
 \inf\{\Ir(\psi):\psi\in \oo{\A^k}\}\geq \Ir(\psi^k)-\frac{\eps}{2},
 \qquad k=1,2,\ldots,N,
 \label{421}
\end{equation}
where, throughout, $\A^k:=\A_{v^k}(\psi^k)$.

Next we define a suitable policy such that the lower bound is asymptotically attained. Let $v=\frac{v_2}{2}\wedge\frac{T}{4}$ and $L=L(v)=\lfloor\frac{T}{v}\rfloor$.
Define $a^\ell=\ell\cdot\frac{T}{L+1}$. Then $[a^\ell, a^{\ell+1})$ forms a disjoint partition of $[0,T)$ with $a^0=0$ and $\sup_i|a^{\ell+1}-a^\ell|<v$. Now consider a sequence $\{\al_n\}$ such
that $\frac{\al_n\sqrt n}{b_n}\to 0$ as $n\to\infty$. Define $H^n=\lfloor \frac{\tilde v}{\al_n}\rfloor$ where $\tilde v=\frac{T}{L+1}$. We can choose $\{\al_n\}$
small enough so that $\inf_nH^n\geq 2$. Define
$$ b^{\ell j}_n=a^\ell+j\frac{\tilde v}{H^n+1}, \quad j=0,\ldots, H^n+1.$$
Hence $\{[b^{\ell j}_n, b_n^{i(j+1)})\}_{j=0}^{j=H^n}$ forms a disjoint partition of $[a^\ell, a^{\ell+1})$ where $b_n^{\ell0}=a^\ell, b_n^{\ell(H^n+1)}=a^{\ell+1}.$
Also $\del_n:=\frac{\tilde v}{H^n+1}<\al_n$ for all $n$. Now we will split each $[b^{\ell j}_n, b_n^{\ell(j+1)}), \ell=0,2,\ldots, L, j=0,\ldots, H^n,$ using some
random intervals. For this we define
\begin{equation}\label{23}
\begin{cases}
\calT^n_i=\int_0^\cdot B^n_i(s)ds,\\
\tilde D^n_i=\tilde S^n_{\mu_i}\circ (\frac{1}{N^n}\calT_i^n),\\
P_n=(\tilde A^n,\tilde D^n).
\end{cases}
\end{equation}
We also denote $\tilde\ell=f(x\cdot\theta)-x$ and
\be\label{26}
F^n_i(a^\ell)=\frac{b_n}{\mu_i\sqrt{n}}\frac{\hat\bzeta_i[P_n](a^{\ell-1})-
\hat\bzeta_i[P_n](a^{\ell-2})}{\tilde v},\quad \ell\geq 2.
\ee
$\theta\cdot\hat\bzeta$ being nondecreasing we have $\sum_i F^n_i(a^\ell)\geq 0$ for all $\ell\geq 2$.
We need to define some more variable before we define the policy. We define random variables $\beta^n_i(\ell), \ell\geq 2,$ as follows:

\[\beta^n_i(\ell)=\left\{
\begin{array}{lll}
(\rho_i-F^n_i(a^\ell))^+\ & \mbox{if} \sum_{i\in\I} (\rho_i-F^n_i(a^\ell))^+\leq 1
\  \mbox{and}\ \|P_n\|^*_{a^{\ell-1}}<M+2,
\\
\rho_i & \mbox{otherwise}.
\end{array}
\right.\]
We also define for $i\in \I$,
\[\gam^n_i=\left\{
\begin{array}{lll}
(\rho_i-\frac{b_n}{\mu^i\sqrt n}\frac{\tilde\ell_i}{\tilde v})^+\ & \mbox{if} \ \sum_{i\in\I} (\rho_i-\frac{b_n}{\mu^i\sqrt n}\frac{\tilde\ell_i}{\tilde v})^+\leq 1
\\
\rho_i & \mbox{otherwise.}
\end{array}
\right.\]
Now we split the intervals $[b^{\ell j}_n, b_n^{\ell(j+1)}), \ell=0,1,\ldots, L, j=0,1,\ldots, H^n,$ using the above variables as follows:
For $\ell=0, j=0,\ldots, H^n$, we define
$$c^{0j}_n(i)=b^{0j}_n+\del_n\sum_{k=1}^i\gam^n_k, \quad i=1,2,\ldots, \bI.$$
We fix the notation as $c^{\ell j}_n(0)=b^{\ell j}_n$ and $c^{\ell j}_n(\bI+1)=b^{\ell(j+1)}$ for all $\ell, j$.
For $\ell=1, j=0,\ldots, H^n$, we define
$$c^{1j}_n(i)=b^{1j}_n+\del_n\sum_{k=1}^i\rho_k, \quad i=1,2,\ldots, \bI.$$
Finally, for $\ell\geq 2, j=0,\ldots, H^n$, we define
$$c^{\ell j}_n(i)=b^{\ell j}_n+\del_n\sum_{k=1}^i\beta^n_k(\ell), \quad i=1,2,\ldots, \bI.$$
It is easy to see that $\{[c^{\ell j}_n(i),c^{\ell j}_n(i+1))\}_{i=0}^{i=\bI}$ forms a partition of $[b^{\ell j}_n, b^{\ell(j+1)}_n)$
for all $\ell=0,\ldots, L,\ j=0,\ldots, H^n$. Now we are ready to define the policy.
Recall from \eqref{3}, \eqref{5} that
\begin{equation}\label{20}
\begin{cases}
D_i^n=S_i\circ (\mu^n_i\calT_i^n),\\
X_i^n=X_i^n(0)+A_i^n-D_i^n.
\end{cases}
\end{equation}
For $i\in\I$, assume $B^n_i$ is given by
\be\label{24}
B^n_i(t)=C^n_i(t)\wedge X^n_i(t), \quad t\in [0, T),
\ee
where
\begin{equation}\label{21}
         C_i^n(t)=\begin{cases}
         \displaystyle
         N^n & \text{if} \ t\in[c^{\ell j}_n(i-1), c^{\ell j}_n(i))\ \text{for some}\ \ell, j
         \\
         0\ & \text{otherwise.}
         \end{cases}
\end{equation}
Now let us argue that $B^n$ is indeed an admissible policy. First we note that $X^n, A^n, D^n$ are piecewise constant processes. Since $\gam^n$ is deterministic,
the policy is well defined on $[0, a^2)$ with RCLL paths. This can be seen by applying induction on the jump times. Since $\beta^n_i(\ell),\ \ell\geq 2,$ is completely determined
by the values of $P_n$ on $[0, a^{\ell-1}]$, $B^n$ is uniquely defined on $[0, T)$. Fix $B^n_i(T)=0$ for all $i\in\I$. It is easy to check that $B^n$ satisfies all
the requirement for being admissible control. Hence $B^n\in \frB^n$ for all $n$. Hence by definition,
\be\label{55}
V^n_X(\tilde X^n(0))\leq J^n_X(\tilde X^n(0), B^n).
\ee

With the policy defined above, we prove the result for the upper bound. In what follows $c_1, c_2,\ldots,$ denote constants independent of $\Del, \eps, \del, \eta, v$ and $n$.

Define $\varphi^k(t)=f(\varphi_\theta[\bar\psi^k](t))$ where $\bar\psi^k=x+yt+\psi^{k,1}(t)-R[\psi^{k,2}](t)$.
Recall from \eqref{018} that $\varphi^k$ is the dynamics corresponding to
$\bar\psi^k$ and $\bzeta[\psi^k]$.
Let $\tilde\La_n=\La_T(\tilde{A}^n) + \La_T(\tilde{S}^n_\mu)$ and denote by
$\Om^n_k$ the event $\{(\tilde{A}^n,\tilde{S}^n_\mu)\in \A^k\}$.
We prove the result in number of steps.
First we show that for a constant $c_1$, for all $n\ge n_0(\eps,v)$,
\begin{equation}\label{41}
\|\tilde X^n\|^*_T\le c_1(1+\tilde\La_n),
\end{equation}
and
\begin{equation}\label{42}
\sup_{[\tilde v,T]}\|\tilde X^n-\varphi^k\|\le c_1\eps,\qquad \text{on } \Om^n_k,\, k=1,2,\ldots,N.
\end{equation}

\spa

\noi{\bf Step 1:}\ From \eqref{311} and the fact $\frac{1}{N^n}\calT^n_i(t)\leq t$, we see that there exists $c_2$ such that
\be\label{59}
\sup_{a^\ell\in[0, t]}\|F^n(a^\ell)\|\leq \frac{b_n}{\sqrt n}\frac{c_2}{\tilde v}(1+\|P_n\|^*_t).
\ee
Since $\rho_i\in(0,1)$ for all $i\in\I$, we note from (\ref{59})
that for all sufficiently large $n$, for any $\ell\geq 2$,
\[
\|P_n\|^*_{a^{\ell-1}}< M+2\quad\text{implies}\quad
\sum_i(\rho_i-F^n_i(a^\ell))^+=\sum_i(\rho_i-F^n_i(a^\ell))\leq 1,
\]
as $\sum_{i}F^n_i(a^\ell)\geq 0$ for all $\ell\geq 2$. Define
$$\hat\tau_n=\min\{\ell\geq 1 :  \|P_n\|^*_{a^{\ell}}\geq M+2\}. $$
First we consider the event $\{\hat\tau_n=1\}$. By definition, $\beta^n_i(\ell)=\rho_i$ for all $\ell\geq 2$ on $\{\hat\tau_n=1\}$.
For all large $n$, given $t\in [b^{\ell j}_n, b^{\ell (j+1)}_n)$, for some $\ell=0,1,\ldots,L,\ j=0,1,\ldots, H^n$, we have on $\{\hat\tau_n=1\},$
\begin{equation}\label{360}
\begin{array}{lll}
&&|\rho_it-\frac{1}{N^n}\int_0^tC^n_i(s)ds|\\
 && \leq |\rho_ia^\ell-\frac{1}{N^n}\int_0^{a^\ell}C^n_i(s)ds|+|\rho_i(b^{\ell j}-a^\ell)-\frac{1}{N^n}\int^{b^{\ell j}}_{a^\ell}C^n_i(s)ds|
 \\ && \, \, \ +\ |\rho_i(t-b^{\ell j}_n)-\frac{1}{N^n}\int_{t-b^{\ell j}_n}^tC^n_i(s)ds|.
 \end{array}
\end{equation}
If $\ell=0$, then the first term on the r.h.s. disappears and the second term is equal to
$$|\rho_i j\del_n-\frac{1}{N^n}j\del_n\gam^n_iN^n|\leq  j\del_n |\frac{b_n}{\mu^i\sqrt n}\frac{\tilde\ell_i}{\tilde v}|\leq |\frac{b_n}{\mu^i\sqrt n}\tilde\ell_i|, $$
where we use the fact that $\frac{j\del_n}{\tilde v}\leq 1$ for $j\leq H^n+1$.
If $\ell\geq 1$, on $\{\hat\tau_n=1\}$, second term on the r.h.s. of \eqref{360} is equal to $0$ and the first term less than equal to $|\frac{b_n}{\mu^i\sqrt n}\tilde\ell_i|$.
Hence using the fact that $C^n_i\leq N^n$ and $|t-b^{\ell j}_n|\leq\del_n$ we have
$$\sup_{t\in[0, T]}\frac{\sqrt n}{b_n}|\rho_it-\frac{1}{N^n}\int_0^tC^n_i(s)ds|\leq c_3,$$
on $\{\hat\tau_n=1\}$ for all $n$ large. Now we consider the event $\{\hat\tau_n>1\}$.  For $t\in [b^{\ell j}_n, b^{\ell (j+1)}_n), t< a^{\hat\tau_n+1},$
for some $\ell=0,1,\ldots,L,\ j=0,1,\ldots, H^n,$ we have

\begin{equation}\label{61}
\begin{array}{lll}
&&|\rho_it-\frac{1}{N^n}\int_0^tC^n_i(s)ds|\\
 && \leq |\rho_ia^\ell-\frac{1}{N^n}\int_0^{a^\ell}C^n_i(s)ds|+|\rho_i(b^{\ell j}-a^\ell)-\frac{1}{N^n}\int^{b^{\ell j}_n}_{a^\ell}C^n_i(s)ds|
 \\ && \, \, \ +\ |\rho_i(t-b^{\ell j}_n)-\frac{1}{N^n}\int_{t-b^{\ell j}_n}^tC^n_i(s)ds|.
 \end{array}
\end{equation}
Now if $\ell\leq 1$, then a similar argument as above holds to bound the r.h.s. of \eqref{61}. So we consider $\ell\geq 2$. Then for all $n$ large
\begin{eqnarray*}
&&|\rho_it-\frac{1}{N^n}\int_0^tC^n_i(s)ds|
\\
&=& |\rho_it-\frac{1}{N^n}\Big(\int_0^{a^1}C^n_i(s)ds+\int_{a^1}^{a^\ell}C^n_i(s)ds
+\int_{a^\ell}^{b^{\ell j}}C^n_i(s)ds+\int_{b^{\ell j}}^tC^n_i(s)ds\Big)|
\\
&\leq &\frac{b_n}{\sqrt n}c_3+ |\rho_i(a^\ell-a^2)-\frac{1}{N^n}N^n\tilde v\sum_{k=2}^{\ell-1}\beta^n_i(k)|
\\
&&\ \ + |\rho_i(b^{\ell j}-a^\ell)-\frac{1}{N^n}N^nj\del_n\beta^n_i(k)|+ 2\del_n.
\\
& \leq &\frac{b_n}{\sqrt n}c_3+ |\tilde v\sum_{k=2}^{\ell-1}F^n_i(a^k)|
+ |j\del_nF^n_i(a^\ell)|+ 2\del_n
\\
& \leq &\frac{b_n}{\sqrt n}c_3+ |\frac{b_n}{\mu_i\sqrt n}(\hat\bzeta_i[P_n](a^{\ell-2})-\hat\bzeta_i[P_n](0))|
+ |j\del_nF^n_i(a^\ell)|+ 2\del_n
\\
&\leq &\frac{b_n}{\sqrt n}c_4(1+\|P_n\|^*_t)+2\del_n.
\end{eqnarray*}
So now we are left with case $t\in [b^{\ell j}, b^{\ell (j+1)}), t\geq a^{\hat\tau_n+1},$
for some $\ell=0,1,\ldots,L,\ j=0,1,\ldots, H^n,$ on the event $\{\hat\tau_n>1\}$. We note that,
\begin{eqnarray*}
&&|\rho_it-\frac{1}{N^n}\int_0^tC^n_i(s)ds|
\\
&& \leq |\rho_ia^{\hat\tau_n+1}-\frac{1}{N^n}\int_0^{a^{\hat\tau_n+1}}C^n_i(s)ds|+|\rho_i(t-a^{\hat\tau_n+1})-\frac{1}{N^n}\int_{a^{\hat\tau_n+1}}^tC^n_i(s)ds|
\\
&& \leq |\rho_ia^{\hat\tau_n+1}-\frac{1}{N^n}\int_0^{a^{\hat\tau_n+1}}C^n_i(s)ds|+|\rho_i(t-b^{\ell j})-\frac{1}{N^n}\int_{b^{\ell j}}^tC^n_i(s)ds|
\\
&& \leq \frac{b_n}{\sqrt n}c_4(1+\|P_n\|^*_t)+4\del_n.
\end{eqnarray*}
Hence combining all the calculations above and making use of the fact that $\frac{\del_n\sqrt n}{b_n}\to 0$ as $n\to\iy$ we have a constant $c_5$
such that
\be\label{62}
\sup_{t\in[0, T]}\frac{\sqrt n}{b_n}|\rho_it-\frac{1}{N^n}\int_0^tC^n_i(s)ds|\leq c_5(1+\tilde\La_n),
\ee
for all $n$ large.

\spa

\noi{\bf Step 2:}\ Now we are ready to prove \eqref{41}. Rewrite \eqref{8} as $(\tilde X^n_i-\eps)^+=\hat Y^n_i+\hat Z^n_i$, where
\begin{eqnarray*}
\hat Y^n_i(t) &=& (\tilde X^n_i-\eps)^-+\tilde{X}^i_n(0) + y_i^nt
 +\tilde{A}_i^n(t)-\tilde{S}_{\mu_i}^n(\frac{1}{N^n}\calT^i_n(t))
 \\
 &&\ \ + \frac{\mu^n_iN^n}{n}\frac{\sqrt n}{b_n}(\rho_it-\frac{1}{N^n}\int_0^tC^n_i(s)ds)
 \\
 \hat Z^n_i(t)&=& \frac{\mu^n_iN^n}{n}\frac{\sqrt n}{b_n}\frac{1}{N^n}\int_0^t(C^n_i(s)-X^n_i(s))^+ds.
\end{eqnarray*}
Now for each $i\in\I$, $(\tilde X^n_i-\eps)^+$ is nonnegative and $\hat Z^n_i$ is nonnegative, nondecreasing.
Since $\frac{N^n}{b_n\sqrt n}\to 0$ as $n\to\iy$, for all $n$ large (depending on $\eps$),
$$(\tilde X^n_i(s)-\eps)^+>0\To X^n_i(s)>\rho_iN^n+\eps\sqrt n b_n\To X^n_i(s)> N^n\To (C^n_i(s)-X^n_i(s))^+=0. $$
Therefore $\int_0^\cdot (\tilde X^n_i(s)-\eps)^+d\hat Z^n_i(s)=0$. Therefore
$((\tilde X^n_i-\eps)^+,\hat Z^i_n)$ is the solution to the
Skorohod problem for data $\hat Y^i_n$. Hence applying Lipschitz property of the Skorohod map and \eqref{62}, we have
\begin{equation}\label{37}
|\hat Z_i^n|^*_T+|(\tilde X_i^n-\eps)^+|^*_T\leq4|\hat Y^i_n|^*_T\le c_6(1+\tilde\La_n),
\end{equation}
for all $n$ large where we used the fact that $(\tilde X^n_i-\eps)^-\leq \eps+\frac{\rho_iN^n}{b_n\sqrt n}<1$. Now \eqref{41} follows from \eqref{37}.
Since
$$ \frac{\mu^n_iN^n}{n}\frac{\sqrt n}{b_n}(\rho_it-\frac{1}{N^n}\int_0^tB^n_i(s)ds)= \frac{\mu^n_iN^n}{n}\frac{\sqrt n}{b_n}(\rho_it-\frac{1}{N^n}\int_0^tC^n_i(s)ds)+\hat Z^n_i,$$
using \eqref{62}, \eqref{37} and convergence of $\frac{\mu^n_iN^n}{n}$, we have
\be\label{38}
\sup_{t\in[0, T]}\frac{\mu^n_iN^n}{n}\frac{\sqrt n}{b_n}|\rho_it-\frac{1}{N^n}\calT^n_i(t)|\leq c_7(1+\tilde\La_n),
\ee
for all $n$ large.

\spa

\noi{\bf Step 3:}\ In particular, \eqref{38} implies that for all $n$ large,
\begin{equation}
 \sup_{t\in[0, T]}|\rho_it-\frac{1}{N^n}\calT_i^n(t)|\leq \frac{\tilde v}{2}
 \label{427}
\end{equation}
holds on the event $\cup_k\Om^n_k$. Therefore using \eqref{457}, \eqref{456} and \eqref{427} one obtains that
for all large $n$ (see (4.37) in \cite{Atar-Bis}),
\be\label{428}
\sup_{t\in[\tilde v, T]}\|\tilde{S}^n_\mu(T_n(t)) -R[\psi^{k,2}](t-\tilde v)\|\leq
\frac{\delta}{2},
\ee
on $\Om^n_k, k=1,2,\ldots, N$.
In the rest of this step, we calculate the difference between $Z^n(t)$ (see \eqref{9}) and $\bzeta[\psi^k](t-\tilde v)$ on the
event $\Om^n_k$. Recall $\hat\tau_n$ defined above. We note that on $\Om^n_k$ one has $\hat\tau_n> L$ for all $n$ large since $\|P_n\|^*_T\leq\|\tilde A^n\|^*+\|\tilde S^n_\mu\|^*<M+2$
by \eqref{457}. Hence for all $n$ large, on $\Om^n_k, k=1,2,\ldots,N$,
$$\gam^n_i=(\rho_i-\frac{b_n}{\mu_i\sqrt n}\frac{\tilde\ell_i}{\tilde v})\quad \text{and}\quad \beta^n_i(\ell)=(\rho_i-F^n_i(a^\ell)),$$
for $\ell=2,3,\ldots, L$. Therefore for any $t\in [b^{0j},b^{0(j+1)}), \ j=1,\ldots,H^n,$ we have
\begin{align}
|\mu_i\frac{\sqrt n}{b_n}(\rho_it-\frac{1}{N^n}\int_0^tC^n_i(s)ds)-\frac{t}{\tilde v}\tilde\ell_i |&=|\mu_i\frac{\sqrt n}{b_n}(\rho_it-
\frac{1}{N^n}\int_0^{b^{0j}}C^n_i(s)ds-\frac{1}{N^n}\int_{b^{0j}}^tC^n_i(s)ds)-\frac{t}{\tilde v}\tilde\ell_i|\notag
\\
&\leq \mu_i\frac{\sqrt n}{b_n}(t-j\del_n)+\frac{t-j\del_n}{\tilde v}\tilde\ell_i+\mu_i\frac{\sqrt n}{b_n}\del_n\notag
\\
&\leq  2\mu_i\frac{\sqrt n}{b_n}\al_n+\frac{1}{H^n+1}\tilde\ell_i.\label{367}
\end{align}
Now for $k=1,2,\ldots,N$, consider
$$
\hat W_{i,k}^n(t):=\mu_i\frac{\sqrt{n}}{b_n}\Big(\rho_it-\frac{1}{N^n}\int_0^tC^n_i(s)ds\Big)
-\bzeta_i[\bar\psi^k](t-\tilde v),\qquad t\in[\tilde v,T],
$$
on the event $\Om^k_n$. We note from (\ref{018}) that $\bzeta[\psi^k](0)=\tilde\ell$. Hence
for $t\in[a^1,a^2)$ (recall $a^\ell=\ell\tilde v$) and all large $n$, we have from (\ref{456}) and (\ref{316}) that
\begin{align}
|\hat W_{i,k}^n(t)|&
\leq |W^n_{i,k}(\tilde v)-\tilde\ell_i|+|\mu^i\frac{\sqrt{n}}{b_n}(\rho_i(t-a^1)-\frac{1}{N^n}\int_{a^1}^tC^i_n(s)ds)|+
|\tilde\ell_i-\bzeta_i[\bar\psi^k](t-v)|\notag
\\
&\leq 4\mu_i\frac{\sqrt n}{b_n}\al_n+\frac{1}{H^n+1}\tilde\ell_i+\eps,\label{368}
\end{align}
where we use \eqref{367} to estimate the first term and a similar estimate to calculate the second.
Now we consider $t\in [a^2, T)$. Let $t\in[b^{\ell j}_n, b^{\ell (j+1)}_n)$ for some $\ell\geq 2, j=0,\ldots,H^n$. The following calculations are of same type
as step 1. We note that for large $n$, on $\Om^n_k$,
\begin{eqnarray*}
&&\mu^i\frac{\sqrt{n}}{b_n}\Big(\rho_ib^{\ell j}_n-\frac{1}{N^n}\int_0^{b^{\ell j}_n}C^n_i(s)ds\Big)
\\
&&= \mu^i\frac{\sqrt{n}}{b_n}\Big(\rho_ia^\ell-\frac{1}{N^n}\int_0^{a^\ell}C^n_i(s)ds\Big)
+\mu^i\frac{\sqrt{n}}{b_n}\Big(\rho_i(a^\ell-b^{\ell j}_n)-\frac{1}{N^n}\int_{a^\ell}^{b^{\ell j}_n}C^n_i(s)ds\Big)
\\
&&= \mu^i\frac{\sqrt{n}}{b_n}\Big(\rho_ia^\ell-\frac{1}{N^n} N^n\tilde v(\gam^n_i+\rho_i+\sum_{k=2}^{\ell-1}\beta^n_i(k)))
+\mu^i\frac{\sqrt{n}}{b_n}\Big(\rho_i(a^\ell-b^{\ell j}_n)-\frac{1}{N^n}N^n j\del_n\beta^n_i(\ell)\Big)
\\
&&=\hat\bzeta_i[P_n]((\ell-2)\tilde v)+\frac{j\del_n}{\tilde v}[\hat\bzeta_i[P_n]((\ell-1)\tilde v)-\hat\bzeta_i[P_n]((\ell-2)\tilde v)].
\end{eqnarray*}
Therefore using \eqref{315}, \eqref{316}, \eqref{456}, \eqref{457} and \eqref{427}, we see that for all $n$ large, on $\Om^n_k, k=1,2,\ldots,N$,
\begin{eqnarray*}
|\hat\bzeta_i[P_n]((\ell-2)\tilde v)-\bzeta_i[\psi^k](t-\tilde v)| &\leq & c_8\eps
\\
|\hat\bzeta_i[P_n]((\ell-1)\tilde v)-\hat\bzeta[P_n]((\ell-2)\tilde v)|&\leq & c_8\eps
\\
|\frac{\sqrt{n}}{b_n}\Big(\rho_it-\frac{1}{N^n}\int_0^{t}C^n_i(s)ds\Big)-
\frac{\sqrt{n}}{b_n}\Big(\rho_ib^{\ell j}_n-\frac{1}{N^n}\int_0^{b^{\ell j}_n}C^n_i(s)ds\Big)| &\leq & 2\frac{\sqrt n}{b_n}\al_n,
\end{eqnarray*}
for some constant $c_8$ independent of $t$ (see (4.40) in \cite{Atar-Bis}). Hence combining all these calculations with \eqref{368} and using the fact that $\frac{\sqrt n}{b_n}\al_n\to 0$,
we have, for all $n$ large and all $k$,
\begin{equation}
\sup_{t\in[\tilde v,T]}\Big|\frac{N^n\mu_i^n}{n}\frac{\sqrt{n}}{b_n}\Big(\rho_it-\frac{1}{N^n}\int_0^tC_i^n(s)ds\Big)
 -\bzeta_i[\bar\psi^k](t-\tilde v)\Big|\le c_9\eps,
 \label{430}
\end{equation}
on $\Om^k_n$.

\spa

\noi{\bf Step 4:}\ Recall $\varphi^k(t)=f(\varphi_\theta[\psi^k](t))$.
The goal of this step is to estimate the difference between $\tilde X_n$ and $\varphi^k$ on $\Om^k_n$.
To this end, let first
\[
\tilde\varphi^k(t)=\left\{\begin{array}{ll}
\ds
x+\frac{t}{\tilde v}\tilde\ell & \mbox{for}\ t\in[0,\tilde v)
\\
f(\varphi_\theta[\psi^k](t-\tilde v)) & \mbox{for}\ t\in[\tilde v,T].
\end{array}
\right.
\]
Also recall from step 2 that $(\tilde X^n_i-\eps)^+$ solves Skorohod problem for the date $\hat Y^n_i$. Since $\Gamma(\tilde\varphi^k)=\tilde\varphi^k$ for all $k=1,\ldots,N$, we have
for large $n$,
\be\label{50}
|(\tilde{X}_i^n-\eps)^+-\tilde\varphi^k_i|^*_T\leq 2|\hat Y^i_n-\tilde\varphi^k_i|^*_T.
\ee
Now for $t\in [0,\tilde v)$,
\begin{eqnarray*}
&&|\hat Y_i^n(t)-\tilde\varphi^k_i(t)|
\\
&&\leq |(\tilde X^n_i-\eps)^-+\tilde{X}^i_n(0) + y_i^nt
 +\tilde{A}_i^n(t)-\tilde{S}_{\mu_i}^n(\frac{1}{N^n}T^i_n(t))
\\
&&\ + \frac{\mu^n_iN^n}{n}\frac{\sqrt n}{b_n}(\rho_it-\frac{1}{N^n}\int_0^tC^n_i(s)ds)-x_i-\frac{t}{\tilde v}\tilde\ell_i|
\\
&&\leq c_{10}\eps,
\end{eqnarray*}
for all $n$ large where we use \eqref{367}, \eqref{456} and \eqref{457}. Similarly, using \eqref{428}, \eqref{430}, for $t\in[\tilde v, T]$,
\begin{eqnarray*}
&&|\hat Y_i^n(t)-\tilde\varphi^k_i(t)|
\\
&&\leq |(\tilde X^n_i-\eps)^-+\tilde{X}^i_n(0) + y_i^nt
 +\tilde{A}_i^n(t)-\tilde{S}_{\mu_i}^n(\frac{1}{N^n}T^i_n(t))-\psi^{k,1}_i(t-\tilde v)+R[\psi^{k,2}]_i(t-\tilde v)
\\
&&\ + \frac{\mu^n_iN^n}{n}\frac{\sqrt n}{b_n}(\rho_it-\frac{1}{N^n}\int_0^tC^n_i(s)ds)-x_i-y_i(t-\tilde v)-\bzeta_i[\psi^k](t-\tilde v)|
\\
&&\leq c_{11}\eps.
\end{eqnarray*}
Therefore from \eqref{50}, we get that for large $n$, on $\Om^n_k, k=1,\ldots,N$,
$$|\tilde X^n_i-\tilde\varphi^k_i|^*_T\leq c_{12}\eps.$$
Thus \eqref{42} follows by comparing $\tilde\varphi^k$ and $\varphi^k$.

Rest of the proof follows by standard argument using \eqref{41} and \eqref{42} (see for example, Step 5 in \cite{Atar-Bis}).\hfill $\Box$

\spa

If $\frac{N^n}{b_n\sqrt n}\to 0$, then $\|\tilde X^n-\tilde Q^n\|^*_T\to 0$ as $n\to\infty$. Hence it is easy to obtain estimates like \eqref{41} and \eqref{42}
for $\tilde Q^n$ with the policy constructed in Theorem \ref{upper1} when $\frac{N^n}{b_n\sqrt n}\to 0$. Thus we have the following theorem:
\begin{theorem}
Let Conditions \ref{moderate}, \ref{mini} and \ref{unbounded} hold and $\lim_{n\to\iy}\frac{N^n}{b_n\sqrt n}=0$. Then
$$\limsup_{n\to\iy} V^n_Q(\tilde{Q}^n(0))\leq V(x).$$
\end{theorem}

\begin{theorem}
Let Conditions \ref{moderate} and \ref{mini} hold. If $h$ and $g$ are bounded, then
$$\limsup_{n\to\iy} V^n_{X}(\tilde{X}^n(0))\leq V(x).$$
\end{theorem}

\noi{\bf Proof:}\ From Theorem \ref{upper1}, we see that we only need to consider the case when $\limsup\frac{N^n}{b_n\sqrt n}>0$.
Instead of introducing a new subsequence, we assume that $\lim\frac{N^n}{b_n\sqrt n}>0$. Hence $\lim\frac{\sqrt n}{b_nN^n}=0$.

Given $\eps\in(0,1)$, we construct an $\eps-$optimal policy. Since $h, g$ are bounded, it is enough
to construct a policy so that \eqref{42} holds for large $n$.

Let $\Del>0$ be given. Define
\be\label{760}
\calQ=\{\psi\in D([0, T], \R^{2\bI})\ : \ \Ir(\psi)\leq \Del\}.
\ee
Hence $\calQ\subset D(M)$ for a suitably chosen $M$.
Using the same argument as in Theorem \ref{upper1}, we have $v\leq\frac{\eps}{2}\wedge\frac{T}{4}$ such that
\begin{equation}
 \osc_{v}(\psi^l_i)<\frac{\delta}{4\sqrt{2\bI}}, \ \mbox{for all}\ \psi=(\psi^1,\psi^2)\in \calQ,\,l=1,2,\,
 i\in\I,\label{7456}
\end{equation}
where $\del\in(0,\eps)$ is chosen according to \eqref{315} and
\begin{equation}
\psi\in \A_{v}(\tilde\psi) \quad \text{ implies } \quad
\|\psi-\tilde\psi\|^*<\frac{\delta}{4},\label{7457}
\end{equation}
for all $\tilde\psi\in\calQ$.
Also we can find finite number of members $\psi^1, \psi^2,\ldots,\psi^N$
of $\calQ$, and positive constants $v^1,\ldots, v^N$ with $v^k<v$, satisfying
$\calQ\subset\cup_k \A^k$, and
\begin{equation}
 \inf\{\Ir(\psi):\psi\in \oo{\A^k}\}\geq \Ir(\bar\psi^k)-\frac{\eps}{2},
 \qquad k=1,2,\ldots,N,
 \label{7421}
\end{equation}
where, throughout, $\A^k:=\A_{v^k}(\bar\psi^k)$. Define
\be\label{726}
F^n_i(t)=\frac{b_n}{\mu_i\sqrt{n}}\frac{\hat\bzeta_i[P_n](j v)-
\hat\bzeta[P_n]((j-1)v)}{v},\quad \mbox{for}\ jv\leq t<(j+1)v,\ j\geq 1,
\ee
where $P_n=(\tilde A^n, \tilde D^n)$ \eqref{23}. Since $\hat\bzeta$ satisfies the causality property, $F^n$ is well defined.
Denote
$$\Theta(a,b)=a\chi_{\R^+}(a)\chi_{[0,1]}(b), \quad a, b\in\R. $$
Recall that $\tilde\ell=f(x\cdot\theta)-x$.
Define
\be
B^n_i(t)=
\begin{cases}
\displaystyle\Theta(\lfloor(\rho_i-\frac{b_n}{\mu_i\sqrt n}\frac{\tilde\ell_i}{v})N^n\rfloor,
\sum_{i\in\I}(\rho_i-\frac{b_n}{\mu_i\sqrt n}\frac{\tilde\ell_i}{v})^+)\chi_{\{\tilde X^n_i(t)>\eps\}}\wedge X^n_i(t) & \text{if}\ t\in[0,v)
\\ \\
\displaystyle\lfloor\rho_iN^n\rfloor\chi_{\{\tilde X^n_i(t)>\eps\}}\wedge X^n_i(t) & \text{if} \ t\in[v, 2v)
\\ \\
\displaystyle\Theta(\lfloor(\rho_i-F^n_i(t-v))N^n\rfloor,\sum_{i\in\I}(\rho_i-F^n_i(t-v))^+)\chi_{\{\tilde X^n_i(t)>\eps\}}\wedge X^n_i(t)
& \text{otherwise}.
\end{cases}
\ee

Using same argument as in Theorem \ref{upper1}, it is easy to see that $B^n$ is an admissible control and hence $B^n\in\frB^n$.
As earlier, define $\varphi^k(t)=f(\varphi_\theta[\psi^k](t))$ where $\bar\psi^k=x+yt+\psi^{k,1}(t)-R[\psi^{k,2}](t)$.
Denote by
$\Om^n_k$ the event $\{(\tilde{A}^n,\tilde{S}^n_\mu)\in \A^k\}$.

\spa

In what follows, $c_1, c_2,\ldots$ denote constants independent of $\Del, \eps, n, v, \del, \eta$.

As earlier (proof of Theorem \ref{upper1}), it is enough to show that
there exists a constant $c_1$ such for all $n\ge n_0(\eps,v)$,

\begin{equation}\label{742}
\sup_{[v,T]}\|\tilde X^n-\varphi^k\|\le c_1\eps,\qquad \text{on } \Om_n^k,\, k=1,2,\ldots,N.
\end{equation}
First we note from \eqref{7457} that
$\|P_n\|^*_T<M+2$ on $\Om^n_k$ for $k=1,2,\ldots, N$. Also
from \eqref{311} and the fact $\frac{1}{N^n}T^n_i(t)\leq t$, we see that there exists $c_2$ such that
\be\label{759}
\sup_{s\in[v, t]}\|F^n(s)\|\leq \frac{b_n}{\sqrt n}\frac{c_2}{v}(1+\|P_n\|^*_t).
\ee
Therefore for all $n$ large, $\rho_i-\frac{b_n}{\mu_i\sqrt n}\frac{\tilde\ell_i}{v}, \rho_i-F^n_i(t-v), i\in\I,$ are
positive and
$$\sum_i(\rho_i-\frac{b_n}{\mu_i\sqrt n}\frac{\tilde\ell_i}{v})\leq 1,\quad \sum_i\rho_i-F^n_i(t-v)\leq 1$$
on $\Om^n_k, k=1,2,\ldots, N$ where we use the fact that $\theta\cdot\tilde\ell=0$ and $\sum_i F^n_i\geq0$.
Again
\begin{eqnarray*}
\tilde X^n_i(t)> \eps\To X^n_i(t)>\rho_i N^n+\eps b_n\sqrt{n}\To \eps\To X^n_i(t)> (\rho_i+\eps\frac{b_n}{\sqrt n}\frac{n}{N^n})N^n.
\end{eqnarray*}
Hence using \eqref{759} and the fact $\lim\frac{n}{N^n}=\iy$, we have on $\Om^n_k, k=1,2,\ldots, N$,
\begin{eqnarray*}
\tilde X^n_i(t)> \eps\To X^n_i(t)> (\rho_i-F^n_i(t-v))N^n, \quad t\geq 2v,
\end{eqnarray*}
for all $n$ large. Similar fact holds in $[0, v)$. Therefore for large $n$, we have
$$B^n_i(t)=\lfloor C^n_i(t)\rfloor\chi_{\{\tilde X^n_i(t)>\eps\}}$$ on $\Om^n_k,\ k=1,2,\ldots,N,$ where
\be\label{790}
C^n_i(t)=
\begin{cases}
\displaystyle(\rho_i-\frac{b_n}{\mu_i\sqrt n}\frac{\tilde\ell_i}{v})N^n & \text{if}\ t\in[0,v)
\\ \\
\displaystyle\rho_iN^n & \text{if} \ t\in[v, 2v)
\\ \\
\displaystyle(\rho_i-F^n_i(t-v))N^n
& \text{otherwise}.
\end{cases}
\ee
Using the same argument as \eqref{666}, we have for large $n$,
\be\label{791}
\displaystyle\sup_{i\in\I}|(\tilde X^n_i)^-|^*_T\leq (c+6)\eps,
\ee
on $\Om^n_k, k=1,2,\ldots,N$ where $c=\sup_n\|y^n\|$. Following the same arguments in (Step 1, \cite{Atar-Bis}) we obtain, for
large $n$,
\be\label{792}
\sup_{[0,T]}\frac{\sqrt n}{b_n}|\rho_i-\frac{1}{N^n}\int_0^tC^n_i(s)ds|\leq c_4,
\ee
on $\Om^n_k$ for some constant $c_4$. We rewrite \eqref{8} as $(\tilde X^n_i-\eps)^+=\hat Y^n_i+\hat Z^n_i$ where
\begin{eqnarray*}
\hat Y^n_i(t) &=& (\tilde X^n_i-\eps)^-+\tilde{X}^i_n(0) + y_i^nt
 +\tilde{A}_i^n(t)-\tilde{S}_{\mu_i}^n(\frac{1}{N^n}T^i_n(t))
 \\
 &&\ \ + \frac{\mu^n_iN^n}{n}\frac{\sqrt n}{b_n}(\rho_it-\frac{1}{N^n}\int_0^tC^n_i(s)ds)
 + \frac{\mu^n_iN^n}{n}\frac{\sqrt n}{b_n}\frac{1}{N^n}\int_0^t(C^n_i(s)-\lfloor C^n_i(s)\rfloor) ds
 \\
 \hat Z^n_i(t)&=& \frac{\mu^n_iN^n}{n}\frac{\sqrt n}{b_n}\frac{1}{N^n}\int_0^t\lfloor C^n_i(s)\rfloor\chi_{\{\tilde X^n_i(s)\leq \eps\}}ds.
\end{eqnarray*}
Since $(\tilde X^n_i-\eps)^+>0\To \chi_{\{\tilde X^n_i(s)\leq \eps\}}=0$, $(\tilde X^n_i-\eps)^+$ solves Skorohod problem
for the data $\hat Y^n_i$. Hence on $\Om^n_k$, for large $n$, $\sup_i|\hat Z^n_i|^*_T\leq c_5$
 using the fact that
$$\sup_{[0, T]}|\frac{\sqrt n}{b_n}\frac{1}{N^n}\int_0^t(C^n_i(s)-\lfloor C^n_i(s)\rfloor) ds|\leq
\frac{\sqrt n}{b_n}\frac{T}{N^n} \to 0 .$$
 Combining
with \eqref{792}, for large $n$,
\be\label{793}
\displaystyle\sup_i\sup_{[0,T]}\frac{\sqrt n}{b_n}|\rho_i-\frac{1}{N^n}\int_0^tB^n_i(s)ds|\leq c_6,
\ee
on $\Om^n_k$ for some constant $c_6$. Now we can use the same arguments as in (Step 3, \cite{Atar-Bis}) to conclude that
for large $n$, on $\Om^n_k$,
\be\label{794}
\sup_{[v, T]}|\frac{N^n\mu^n_i}{n}\frac{\sqrt n}{b_n}(\rho_it-\frac{1}{N^n}\int_0^tB^n_i(s)ds)-\bzeta_i[\bar\psi^k](t-v)|
\leq c_7\eps,
\ee
for some constant $c_7$. We define $\tilde\varphi^k$ as in Step 4 above replacing $\tilde v$ by $v$. Using the Lipschitz
property of the Skorohod map, we obtain on $\Om^n_k$, for large $n$,
$$|(\tilde X^n_i)^+-\tilde\varphi^k_i|^*_T\leq2|\hat Y^n_i-\tilde\varphi^k_i|^*_T\leq c_8\eps,$$
where the last estimate is obtained using the same argument as Step 4 above. Combining with \eqref{791}, we have on
$\Om^n_k, k=1,2, \ldots, N$, $|\tilde X^n_i-\tilde\varphi^k_i|^*_T\leq c_9\eps$ for some constant $c_9$ and $n$ large.
Hence we obtain \eqref{742} comparing $\tilde\varphi$ and $\varphi$.\hfill $\Box$

\begin{remark}\label{end-rem}
One can weaken the Condition \ref{mini} by assuming the existence of two continuous minimizer $f_h, f_g,$ corresponding to $h$ and $g$,
respectively. It is possible to show that there exists a minimizing strategy $\tilde\bzeta$ in this setting which is equal to $\bzeta$ in $[0, T)$ and takes a jump at time $T$ (\cite[Remark 3.1]{Atar-Bis}). $\tilde\bzeta$ will also have similar regularity properties as $\bzeta$. Therefore
the proof of the lower bound will be very similar to what we have presented here. One needs to modify the policy for the upper bound in the interval $[T-\tilde v, T] (or\ [T-v, T])$ to incorporate the jump at time $T$. This can be done in the same manner as we treat the jump at time $0$. However we do not add all the details here for simplicity. 
\end{remark}

\subsection{Linear cost and asymptotic optimality}\label{sec3.3}
In this section, we provide a simple policy based on priority that is asymptotically optimal. We assume that $h$ and $g$ have following forms
\[
h(x)=\sum_{i=1}^\bI c_ix_i,\qquad g(x)=\sum_{i=1}^\bI d_ix_i,
\]
where $c_i$ and $d_i$ are nonnegative constants, and, in addition,
\[
c_1\mu_1\geq c_2\mu_2\geq\cdots\geq c_\bI\mu_\bI \quad \text{and} \quad
d_1\mu_1\geq d_2\mu_2\geq\cdots\geq d_\bI\mu_\bI.
\]
We consider the $c\mu$-rule that prioritizes according to the ordering of class labels, with highest priority to class 1. Define
\be\label{555}
B^n_1=X^n_1\wedge N^n,\ B^n_2=X^n_2\wedge(N^n-B^n_1),\ldots,\ B^n_\bI=X^n_\bI\wedge(N^n-\sum_{i<\bI}B^n_i).
\ee
This policy is in the spirit of the priority policy considered by Cox and Smith \cite[Chapter III]{cox-smith} for the linear cost
in a single server queuing network.
It is easy to see that the above policy is consistent with \eqref{3}-\eqref{5} and $B^n\in\frB^n$.
Proof of the following theorem follows using the same argument from Theorem 5.1 in \cite{Atar-Bis}.
\begin{theorem}
Assume Conditions \ref{moderate}, \ref{unbounded} hold and $\lim\frac{N^n}{b_n\sqrt n}=0$.
Then, under the priority policy $\{B^n\}$ in \eqref{555},
$$
\lim_{n\to\infty}J^n_Q(\tilde{Q}^n(0),B^n)=\lim_{n\to\infty}J^n_X(\tilde{X}^n(0),B^n)=V(x).
$$
\end{theorem}

\spa

\noi{\bf Acknowledgement:}\ The author is grateful to Prof. Rami Atar for valuable discussions.


\bibliographystyle{is-abbrv}

\bibliography{refs}

\end{document}